\documentclass[11pt]{article}
\usepackage{amsfonts}
\usepackage{latexsym,amssymb,amsmath,graphics,cite,arydshln}
\usepackage{tikz}
\usepackage{fancyhdr}
\usepackage{lipsum}
\usepackage{lineno}
\usepackage{color}
\usepackage{diagbox}
\usepackage{colortbl}
\usepackage{multirow}

\usepackage[margin=1in]{geometry}
\usepackage{amsmath,amssymb,amsthm,mathtools}
\usepackage{microtype}
\usepackage{enumitem}
\usepackage{booktabs}
\usepackage{array}
\usepackage[hidelinks]{hyperref}
\usepackage[nameinlink,capitalise,noabbrev]{cleveref}
\usepackage{lmodern}

\begin{document}

\newtheorem{theorem}{Theorem}[section]
\newtheorem{lemma}[theorem]{Lemma}
\newtheorem{corollary}[theorem]{Corollary}
\newtheorem{remark}[theorem]{Remark}
\newtheorem{example}[theorem]{Example}
\newtheorem{definition}[theorem]{Definition}
\newtheorem{construction}[theorem]{Construction}
\newtheorem{fact}[theorem]{Fact}
\newtheorem{proposition}[theorem]{Proposition}
\newtheorem{conjecture}[theorem]{Conjecture}
\newtheorem{claim}[theorem]{Claim}

\newcommand{\Z}{\mathbb Z}
\newcommand{\N}{\mathbb N}
\newcommand{\E}{\mathbb E}
\newcommand{\Prb}{\mathbb P}
\newcommand{\dist}{\operatorname{dist}}
\newcommand{\cl}{\operatorname{cl}}
\newcommand{\supp}{\operatorname{supp}}
\newcommand{\dd}{\,\mathrm d}
\newcommand{\K}{K}
\newcommand{\Awt}{\mathcal A_{w,t}}
\newcommand{\Shadow}{\partial_2}
\newcommand{\tp}{\tau}
\newcommand{\pospart}[1]{\left(#1\right)_{+}}

\newenvironment{poc}{\begin{proof}[Proof of the claim]\renewcommand*{\qedsymbol}{$\blacksquare$}}{\end{proof}}

\begin{center}
{\Large\bf Optimal bounds for $(n,t,w,w-2)$-constant-power error-correcting cooling codes}
\renewcommand{\thefootnote}{}
\footnotetext{This work is supported by the Natural Science Foundation of the Jiangsu Higher Education Institutions of China under Grant 25KJB110016 (S. Liu), the National Natural Science Foundation of China under Grants $12271023$ (S. Yu), $12601639$ (M. Zhang) and $12671391$ (L. Ji).}

\vskip12pt
Shuangqing Liu$^1$, Qianfei Zhang$^1$, Shikang Yu$^2$, Menglong Zhang$^3$, Lijun Ji$^4$  \\[2ex]
 {\footnotesize $^1$Department of Mathematics, Suzhou University of Science and Technology, Suzhou 215009, P. R. China}\\
 {\footnotesize $^2$School of Mathematics and Statistics, Beijing Jiaotong University, Beijing, 100044, P.R. China}\\
 {\footnotesize$^3$Institute of Mathematics and Interdisciplinary Sciences, Xidian University, Xi'an 710126, P.R. China}\\
 {\footnotesize $^4$Department of Mathematics, Soochow University, Suzhou 215006, P. R. China} \\
{\footnotesize
shuangqingliu@usts.edu.cn, 18014877519@163.com, healthyu@bjtu.edu.cn, zhangmenglong@xidian.edu.cn, jilijun@suda.edu.cn}
\vskip12pt
\end{center}

\vskip12pt

\noindent {\bf Abstract:}
High temperatures have dramatic negative effects on interconnect performance. In a bus, whenever a wire changes its state from ``0'' to ``1'' or from ``1'' to ``0'', Joule heating causes the temperature to rise. Constant-power error-correcting cooling (CPECC) codes, introduced in [IEEE Trans. Inf. Theory, 66 (2020), 4804--4818], are coding schemes that can simultaneously control the peak temperature and the average power consumption of on-chip buses while providing error-correction capability for transmitted information. This paper establishes a new upper bound for $(n,t,w,w-2)$-CPECC codes for $t<w-3$ by some graph-theoretic techniques and the probabilistic method, and constructs several new families of optimal CPECC codes using the asymptotic existence result for GDDs. For $t\geq w-3$, this paper also provides a new upper bound for $(n,t,w,w-2)$-CPECC codes by some hypergraph-theoretic techniques, and several new families of CPECC codes attaining the resulting upper bounds using RBIBDs. The previous upper bound for $(n,t,w,w-2)$-CPECC codes in [arXiv:2608.08742] is a special case of our upper bound.


\section{Introduction}

Power and heat dissipation limits have emerged as first-order design constraints for chips, whether targeted for battery-powered devices or high-end systems. Numerous encoding techniques have been proposed in the literature \cite{bmmss,ccl,po,sm} in order to reduce the overall power consumption of both on-chip and off-chip buses, but these techniques do not directly address peak temperature minimization. Power-aware design alone is insufficient to address thermal challenges, since it does not directly target the spatial and temporal behavior of the operating environment.

Chee et al. \cite{1} introduced a new class of codes, called cooling codes, to directly control the peak temperature of a bus by cooling its hottest wires. This is achieved by avoiding state transitions on the hottest wires for as long as necessary until their temperature drops. Cooling codes are based on differential encoding. Specifically, if the current state of the wires is $(s_1,s_2,\ldots,s_n)$, i.e., wire $i$ is in state $s_i$, and we want to transmit the binary vector $(x_1,x_2,\ldots,x_n)$, we set the wires to the state $(s'_1,s'_2,\ldots,s'_n)$, where $s'_i$ is equal to $s_i+x_i$ modulo $2$ for $1\leq i\leq n$. Therefore, there is a state transition on wire $i$ if and only if $x_i=1$. Chee et al. \cite{1} introduced low-power error-correcting cooling (LPECC) codes, which simultaneously control the peak temperature and the average power consumption and provide error correction for the transmitted information. In addition, if an LPECC code can provide the same average power consumption in every transmission,  then it is called a constant-power error-correcting cooling (CPECC) code \cite{cekvw}.

\begin{center}
{\small \centerline{\footnotesize TABLE I:  KNOWN AND NEW UPPER/LOWER BOUNDS OF CPECC CODES}
\vspace{0.08cm}
\begin{tabular}{|c|c|c|}
\hline $C(n,t,w,e)$/$C'(n,t,w,e)$ & Constraints & Reference \\ \hline
 $C'(n,t,w,w-3)\le \frac{{n\choose 3}}{{w+t\choose3}}$ & $n\ge w+t$, $w\geq w_0(t)$ & \cite{t-h}\\ \hline
 $\lim\limits_{n\rightarrow \infty}C'(n,t,w,w-3)= \frac{{n\choose 3}}{{w+t\choose3}}$ & $n\ge w+t$, $w\geq w_0(t)$ & \cite{t-h}\\ \hline
$C'(qw,q-1,w,e)\geq q^{w-e-1}$ & $q$ is a prime power & \cite{cekvw} \\ \hline
\multirow{3}{*}{ $\lim\limits_{n\rightarrow \infty} C'(n,t,w,e)= \frac{\binom{n}{w-e}}{\binom{w+t}{w-e}}$} & $e\leq w-1$, & \multirow{3}{*}{\cite{z-z}} \\
&  $2\binom{w}{w-e}\geq \binom{w-t-1}{w-e}+\binom{w+t}{w-e}$,  &  \\
&  $\binom{w+t-i}{w-e-i}|\binom{n-i}{w-e-i}$, $0\leq i\leq w-e-1$ &  \\  \hline
\multirow{2}{*}{$C'(n,t,4,2)= \frac{n(n-1)}{10(t+1)}$}& $t=2\ell-1$, $5\ell|n$, $n\equiv  5 \pmod {20}$,
& \multirow{2}{*}{ {\cite{lyj}}} \\
& $n\neq45, 345, 465, 645$ & \\  \hline
{$C'(n,1,w,w-2)\leq h_1 $}& $w\geq 4$ & \cite{ltz}
 \\  \hline
 \multirow{3}{*}{$C'(n,1,w,w-2)= h_2$}& $n\geq e^{{(w+1)}^{{(w+1)}^{2}}}-1$,
& \multirow{3}{*}{ \cite{ltz}} \\
& $w\geq4$, $n\equiv 0\pmod {w}$,& \\
& $n(n+1)\equiv 0\pmod {w(w+1)}$& \\  \hline
{$C'(n,t,w,w-2)\leq  \frac{n(n-1)}{2(2w-3)(t+1)} $}& None & \cite{ltz}
 \\  \hline
 \multirow{4}{*}{$C'(n,t,w,w-2)= \frac{n(n-1)}{2(2w-3)(t+1)}$} &$n> e^{e^{(2w-3)^{12(2w-3)^2}}}$,   &\multirow{4}{*}{\cite{ltz}} \\
& $t=l(w-2)-1$,& \\
& $n\equiv 1\pmod {(2w-4)}$, & \\
& $n\equiv 0\pmod {l(2w-3)}$ & \\ \hline
\multirow{2}{*}{$C'(n,t,w,w-2)= \left\lfloor\frac{n(n-1)}{2(2w-3)(t+1)}\right\rfloor $} & $n$ is enough large,
& \multirow{2}{*}{ \cite{ltz}} \\
& $(i$) or $(ii)$ or $(iii)$ in \cite{ltz}& \\  \hline
$C'(n,t,w,w-2)\leq h_3$ & $n\geq t+w$, $t<w-3$ & Theorem \ref{thm:ub_C'(n,1,w,w-2)} \\ \hline
\multirow{4}{*}{$C'(n,t,w,w-2)= h_4$} & $n$ is large enough, $w>t+3$, $t\equiv1\pmod2$, &  \\
&  $n-w-\frac{t-1}{2}\equiv0\pmod{w+t-1}$,  &Theorem \ref{cor-2} \\
&  $n(n-w-\frac{t-1}{2})\equiv0\pmod{(w+t)(w+t-1)}$ & \\
& and $n\equiv0\pmod{w+\frac{t-1}{2}}$ & \\ \hline
$C'(n,t,w,w-2)\leq\frac{{n\choose 2}}{(2w-3)(t+1)+z_{w-1}(t+1)}$ & $n\geq t+w$, $w\geq4$, $t\geq w-3$ & Theorem \ref{thm:upb_(n,t,w,w-2)-CPECC} \\ \hline
\multirow{2}{*}{$C'(n,t,w,w-2)\leq\frac{n(n-1)}{2e_{w,t}}$} & $n\geq t+w$, $w\geq4$,  & Theorem \ref{thm:upb_(n,t,w,w-2)-CPECC_new}\\
&  $t\geq \max\{w-3,(w-3)\lfloor\frac{w-2}{3}\rfloor-1\}$ & \\\hline
\multirow{4}{*}{$C'(n,t,w,w-2)=\frac{n(n-1)}{2((2w-3)(t+1)+l)}$} & $w\geq4$, $t\geq\max\{(w-3)\lfloor\frac{w-2}{3}\rfloor-1,w-3\}$, &  \\
&  $t+1=(w-3)l$ for $3l\leq w-1$,  &Theorem \ref{cor:rbibd-2} \\
&   $n>e^{e^{(2w-4)^{12(2w-4)^2}}}$, $n\equiv1\pmod{2w-5}$,  & \\
& and $n\equiv0\pmod{l(2w-4)}$ & \\ \hline
\end{tabular}}\end{center}
\smallskip
where $h_1=\left\lfloor\frac{n}{w+1}\left\lfloor\frac{n-1}{w}\right\rfloor\right\rfloor+\Big\lfloor\frac{n(n-1)-
	\left\lfloor\frac{n}{w+1}\left\lfloor\frac{n-1}{w}\right\rfloor\right\rfloor(w+1)w}{2w(w-1)}\Big\rfloor$, $h_2=\left\lfloor\frac{n}{w+1}\left\lfloor\frac{n-1}{w}\right\rfloor\right\rfloor+\left\lfloor\frac{n}{2w}\right\rfloor $, $h_3=\left\lfloor\frac{n}{w+t}\Big\lfloor\frac{n-1}{w+t-1}\Big\rfloor\right\rfloor+
\left\lfloor\frac{n(n-1)-\left\lfloor\frac{n}{w+t}\big\lfloor\frac{n-1}{w+t-1}\big\rfloor\right\rfloor(w+t)(w+t-1)}{2\left(\binom{w+\lfloor(t-1)/2\rfloor}{2}+\binom{w+\lceil(t-1)/2\rceil}{2}\right)}\right\rfloor$,
$h_4=\frac{n(n-w-\frac{t-1}{2})}{(w+t)(w+t-1)}+\Big\lfloor\frac{n}{2w+t-1}\Big\rfloor.$

In this paper, we are particularly interested in CPECC codes. An $(n,t,w,e)$-CPECC code \cite{cekvw} is a coding scheme for communication over a bus consisting of $n$ wires with the following three features:
\begin{itemize}
\item[\textbf{A}($t$):] every transmission causes no state transitions on the $t$ hottest wires;
\item[\textbf{B}($w$):] the total number of state transitions on all the wires is exactly $w$ in every transmission;
\item[\textbf{C}($e$):] up to $e$ transmission errors (0 received as 1, or 1 received as 0) on the $n$ wires can be corrected.
\end{itemize}
If Property $\textbf{B}(w)$ is replaced with $\textbf{B}'(w)$: the total number of state transitions on all the wires is at most $w$ in every transmission, then the CPECC code is called an $(n,t,w,e)$-LPECC code. Note that if the coding scheme only satisfies Property $\textbf{A}(t)$, the code is an $(n,t)$-cooling code.

Chee et al. \cite{1} introduced three methods to construct LPECC codes using complete hypergraphs, dual codes and sunflowers. The first and the last authors \cite{lj} presented some upper bounds and several families of optimal LPECC codes using combinatorial configurations. Chee et al. \cite{cekvw} also constructed CPECC codes using some special linear codes, for example, Reed-Solomon codes. Zhao and Zhang \cite{z-z} also gave some asymptotic results for CPECC/LPECC codes. The first author et al. \cite{lyj,ltz} established new upper bounds for $(n,t,w,w-2)$-CPECC/LPECC codes, as well as some new optimal families. Recently, Tong and Hu in \cite{t-h} showed some good results on $(n,t,w,w-3)$-CPECC/LPECC codes. This paper establishes a new upper bound for $(n,t,w,w-2)$-CPECC codes for $t<w-3$ by some graph-theoretic techniques and the probabilistic
method, and constructs several new families of optimal CPECC codes using the asymptotic existence result for GDDs. For $t\geq w-3$, this paper also provides a new upper bound for $(n,t,w,w-2)$-CPECC codes by some hypergraph-theoretic techniques, and several new families of CPECC codes attaining the resulting upper bounds using RBIBDs.
The known and the new ones are listed in Table I.

The rest of this paper is organized as follows. In Section 2, we introduce the CPECC codes, some graph-theoretic concepts and some related combinatorial configurations.
In Section 3, we first show a tight upper bound for $(n,t,w,w-2)$-CPECC codes for $t<w-3$ (see Theorem \ref{thm:ub_C'(n,1,w,w-2)}) by some graph-theoretic techniques and the probabilistic method. Applying the asymptotic existence result for GDDs, we also prove that when $n$ is sufficiently large, this upper bound is attainable (see Theorems \ref{cor-2} and \ref{cor-gdd-1}). In Section 4, we present an upper bound for $(n,t,w,w-2)$-CPECC codes for $t\geq w-3$ (see Theorem \ref{thm:upb_(n,t,w,w-2)-CPECC}) by some hypergraph-theoretic techniques and a new shadow bound for uniform hypergraphs. By computing $z_{w-1}(t+1)$ in Theorem \ref{thm:upb_(n,t,w,w-2)-CPECC}, we present a combinatorial formulation of the upper bound for $(n,t,w,w-2)$-CPECC codes (see Theorem \ref{thm:upb_(n,t,w,w-2)-CPECC_new}). This result shows the upper bound in \cite{ltz} (see Theorem \ref{thm:ubj_(n,t,w,w-2)-CPECC} (1)) is a special case of our upper bound (see Remark \ref{remark-1}). We also construct some optimal families of CPECC codes by using resolvable BIBDs (see Theorems \ref{new-22} and \ref{cor:rbibd-2}). In Section 5, we conclude this paper.

\section{Preliminaries}

Given a positive integer $n$, the set $\{1,2,\ldots,n\}$ is abbreviated as $[n]$. For a vector $\boldsymbol v=(v_1,v_2,\ldots,v_n)\in\mathbb{F}^{n}_2$, the {\em support} of $\boldsymbol v$ is defined as $supp(\boldsymbol v)\triangleq\{i\in[n]:v_i\neq0\}$ and the {\em Hamming weight} of a vector $\boldsymbol v$ is defined as wt$(\boldsymbol v)\triangleq|supp(\boldsymbol v)|$. The {\em Hamming distance} between $\boldsymbol v=(v_1,v_2,\ldots,v_n)\in\mathbb{F}^{n}_2$ and $\boldsymbol u=(u_1,u_2,\ldots,u_n)\in\mathbb{F}^{n}_2$, denoted by $d(\boldsymbol v,\boldsymbol u)$, is the number of positions where $\boldsymbol v$ and $\boldsymbol u$ differ, that is, $d(\boldsymbol v,\boldsymbol u)=|\{1\leq i\leq n:v_i\neq u_i\}|$.

\begin{definition}
For positive integers $n,t,w,e$ with $w>e$ and $n\geq t+w$, let ${\cal C}_1,{\cal C}_2,\ldots,{\cal C}_M$ be disjoint subsets of $\mathbb F_2^n$.
A set $\mathcal{C}={\cal C}_1\cup{\cal C}_2\cup\dots\cup{\cal C}_M$ is called an $(n,t,w,e)$\textbf{-CPECC code} of size $M$ if $\mathcal{C}$ satisfies the following properties:
\begin{itemize}
   \item[\textbf{A}$(t)$:] for any subset $S\subset[n]$ with $|S|=t$ and any $i\in[M]$, there exists a codeword $\boldsymbol v\in{\cal C}_i$ such that $supp(\boldsymbol v)\cap S=\varnothing$;

   \item[\textbf{B}$(w)$:] {\em wt}$(\boldsymbol v)= w$ for any $\boldsymbol v\in\mathcal C$; and

   \item[\textbf{C}$(e)$:] $d(\boldsymbol v,\boldsymbol u)\geq2e+1$ for any $1\leq i<j\leq M$, $\boldsymbol v\in\mathcal C_i$ and $\boldsymbol u\in\mathcal C_j$.
\end{itemize}
Denote by $C'(n,t,w,e)$ the maximum size of an $(n,t,w,e)$-CPECC code.
\end{definition}

As a binary vector $\boldsymbol v\in\mathbb{F}^{n}_2$ naturally corresponds to a subset of the $n$-element set $[n]$, we give the following combinatorial characterizations of CPECC codes. For an $(n,t,w,e)$-CPECC code $\mathcal{C}={\cal C}_1\cup{\cal C}_2\cup\ldots\cup{\cal C}_M$, let ${\cal P}_i=\{supp(\boldsymbol{v})\colon \boldsymbol{v}\in{\cal C}_i\}$ for each $i\in[M]$ and ${\cal B}=\{supp(\boldsymbol{v})\colon \boldsymbol{v}\in{\cal C}\}$. Note that ${\cal P}_1,{\cal P}_2,\dots,{\cal P}_M$ form a partition of ${\cal B}$. Such a set ${\cal B}$ together with a partition into $M$ parts ${\cal P}_1,\ldots,{\cal P}_M$ is a set-theoretic representation of an $(n,t,w,e)$-CPECC code. Formally, a set-theoretic representation of an $(n,t,w,e)$-CPECC code of size $M$ consists of a set ${\cal B}$ of subsets of the $n$-element set $[n]$ together with a partition of ${\cal B}$ into ${\cal P}_1,\ldots,{\cal P}_M$ such that the following properties hold:
\begin{itemize}
	\item[\textbf{(a)}] (Property $\textbf{A}(t)$)~~for any $T\subset[n]$ with $|T|=t$ and any $i\in[M]$, there exists $B\in\mathcal P_i$ satisfying
      $B\cap T=\varnothing$;
	
	\item[\textbf{(b)}] (Property $\textbf{B}(w)$)~~$|B|=w$ for any $B\in{\cal B}$;

	\item[\textbf{(c)}] (Property $\textbf{C}(e)$)~~$|B_i\Delta B_j|\geq2e+1$ for any $1\leq i<j\leq M$, $B_i\in\mathcal P_i$ and $B_j\in\mathcal P_j$,
      where $B_i\Delta B_j=(B_i\cup B_j)\setminus(B_i\cap B_j)$.
\end{itemize}

Next, we introduce some combinatorial configurations that will be used in the following sections.
A {\em group divisible design} ($k$-GDD) is a triple $(V,\mathcal{G},\mathcal{B})$, where $V$ is an $n$-set of points, $\mathcal{G}$ is a set of subsets of $V$ (called {\em groups}) that form a partition of $V$, and $\mathcal{B}$ is a collection of $k$-subsets of $V$ (called {\em blocks}) such that each pair of distinct points from $V$ occurs either in some group or in exactly one block, but not both. A $k$-GDD of type $g_1^{a_1}g_2^{a_2}\cdots g_u^{a_u}$ is a $k$-GDD in which there are $a_i$ groups of size $g_i$, $i=1,2,\ldots,u$. In particular, a $k$-GDD of type $1^n$ is called a {\em balanced incomplete block design} ($(n,k,1)$-BIBD); in this case, we usually omit the group set. $k$-GDDs are widely used to construct other combinatorial configurations, and their constructions have received much attention; see \cite{Ge2007}. The following theorem gives an asymptotic existence result for $k$-GDDs.

\begin{theorem}[\cite{ljq}]\label{thm:asy_dgg}
For fixed integers $k\geq 2$ and $g\geq 1$, there exists an integer $u_0$ such that, for every integer $u\geq u_0$ satisfying $g^2u(u-1)\equiv 0\pmod{k(k-1)}$ and $g(u-1)\equiv 0\pmod{k-1}$, there exists a $k$-GDD of type $g^u$.
\end{theorem}

An $(n,k,1)$-BIBD $(V,\mathcal{B})$ is called {\em resolvable} if its block set $\mathcal{B}$ admits a partition into parallel classes, each being a partition of the point set $V$. A resolvable $(n,k,1)$-BIBD is abbreviated as an $(n,k,1)$-RBIBD.

The following theorem gives an asymptotic existence result for resolvable balanced incomplete block designs.

\begin{theorem}[\cite{chang}]\label{thm:asy_RBIBD}
For any integer $k\geq 3$, there exists an $(n,k,1)$-RBIBD for every integer $n>e^{e^{k^{12k^2}}}$ satisfying $n\equiv 0\pmod{k}$ and $n\equiv 1\pmod{k-1}$.
\end{theorem}

An {\em $(n,k,1)$-packing} is a pair $(V,\mathcal{B})$, where $V$ is an $n$-set of elements (called {\em points}) and $\mathcal{B}$ is a collection of $k$-subsets of $V$ (called {\em blocks}) such that each pair of distinct points from $V$ occurs in at most one block. The following lemma gives the first Johnson bound for an $(n,k,1)$-packing.

\begin{lemma}[\cite{cd}]\label{lem:J_bound}
Each $(n,k,1)$-packing with $b$ blocks satisfies $b\leq \big\lfloor\frac{n}{k}\big\lfloor\frac{n-1}{k-1}\big\rfloor\big\rfloor$.
\end{lemma}

For a graph $G$, let $V(G)$ and $E(G)$ denote the vertex set and edge set of $G$, respectively, and let $v(G)=|V(G)|$ and $e(G)=|E(G)|$. Let $d_G(v)$ be the {\em degree} of the vertex $v$ in $G$, and let $\delta(G)=\min_{v\in V(G)}d_G(v)$ (resp. $\Delta(G)=\max_{v\in V(G)}d_G(v)$) be the minimum (resp. maximum) degree of $G$.

The {\em complement} $\overline G$ of a graph $G$ has vertex set $V(G)$, and two distinct vertices are adjacent in $\overline G$ if and only if they are not adjacent in $G$. For $S\subseteq V(G)$, let $G[S]$ be the induced subgraph of $G$ on $S$ and $G-S=G[V(G)\setminus S]$. Let $K_r$ be the {\em clique} on $r$ vertices. A set $T\subseteq V(G)$ is a {\em $K_r$-transversal} if it meets every copy of $K_r$ in $G$. A graph is {\em $K_r$-free} if it contains no subgraph isomorphic to $K_r$.

An {\em independent set} is a set of pairwise nonadjacent vertices. For a graph $G$, the {\em clique number} $\omega(G)=\max\{v(Q):Q\ \text{is a clique of } G\}$ and the {\em independence number} $\alpha(G)=\max\{|I|:I\ \text{ is an independent set of } G\}$. Directly from the definition of the complement, we have
\begin{align}\label{eq:complement-alpha}
\alpha(\overline G)=\omega(G).
\end{align}

A {\em vertex cover} of $G$ is a set $C\subseteq V(G)$ that contains at least one endpoint of every edge of $G$. Its minimum possible size is the {\em vertex-cover number}, denoted by $\tau_v(G)$. A set $C$ is a vertex cover if and only if $V(G)\setminus C$ is an independent set. Therefore, for every finite graph $G$, the following holds
\begin{align}\label{eq:cover-alpha}
\tau_v(G)=|V(G)|-\alpha(G).
\end{align}
Also, if $J$ is a subgraph of $G$ and $C$ is a vertex cover of $G$, then $C\cap V(J)$ covers all edges of $J$. Hence
\begin{align}\label{eq:cover-monotone}
\tau_v(J)\leq\tau_v(G).
\end{align}

For any $e\in E(G)$, $G-e$ denotes the graph obtained from $G$ by deleting the edge $e$. An edge $e$ of a graph $G$ is \emph{critical for the vertex-cover number} if $\tau_v(G-e)=\tau_v(G)-1$. A graph is called \emph{$\tau_v$-critical} if every edge is critical for the vertex-cover number. Note that deleting one edge can lower the vertex-cover number by at most one, so the condition $\tau_v(G-e)<\tau_v(G)$ is equivalent to the above equality. The following classical result is due to Erd\H{o}s and Gallai \cite{erdos-gallai}.

\begin{theorem}[Erd\H{o}s--Gallai, \cite{erdos-gallai}]\label{thm:EG}
Let $G$ be a graph with no isolated vertices and $\tau_v(G)\geq1$. If $G$ is $\tau_v$-critical, then $|V(G)|\leq2\tau_v(G)$.
\end{theorem}

Let $k\geq1$ be an integer. A graph $G$ is \emph{$k$-colorable} if there is a map $c:V(G)\longrightarrow\{1,\ldots,k\}$ such that $c(u)\neq c(v)$ whenever $uv\in E(G)$. Such a map is called a {\em proper $k$-coloring}. The {\em chromatic number} $\chi(G)$ of $G$ is the least possible $k$. Note that $k$-colorable means colorable with at most $k$ colors, not necessarily with exactly $k$ colors. If $G$ is $k$-colorable, then $V(G)$ is the union of at most $k$ maximal independent sets. Since all vertices of a clique must receive different colors, every $k$-colorable graph is $K_{k+1}$-free.

\begin{theorem}[Brooks' theorem, \cite{Brooks1941}]\label{thm:brooks}
Let $G$ be a connected graph with maximum degree $\Delta(G)=\Delta\geq3$. If $G$ is neither a clique nor an odd cycle, then $\chi(G)\leq\Delta.$
\end{theorem}

A graph $G$ is {\em $q$-degenerate} if every nonempty induced subgraph of $G$ has a vertex of degree at most $q$ within that induced subgraph. We conclude this section with the following simple yet useful result.
\begin{lemma}\label{lem:degenerate-colour}
Every $q$-degenerate graph is $(q+1)$-colourable.
\end{lemma}

\begin{proof}
Repeatedly remove a vertex whose degree in the current induced subgraph is at most $q$. This produces an ordering $v_1,\dots,v_n$ such that each $v_i$ has at most $q$ neighbours among $v_{i+1},\dots,v_n$. Colour the vertices in the reverse order $v_n,v_{n-1},\dots,v_1$. When $v_i$ is coloured, at most $q$ of its already coloured neighbours forbid colours, so at least one of the $q+1$ colours remains available. This gives a proper $(q+1)$-colouring.
\end{proof}

\section{Optimal $(n,t,w,w-2)$-CPECC codes for $t<w-3$}\label{sec:t<w-3}

In this section, we establish a new upper bound for $(n,t,w,w-2)$-CPECC codes for $t<w-3$, which generalizes the upper bound for $(n,1,w,w-2)$-CPECC codes in \cite{ltz}, and obtain new families of optimal $(n,t,w,w-2)$-CPECC codes.

Let $\mathcal B=\mathcal P_1\cup\dots\cup\mathcal P_M$ be the set-theoretic representation of an $(n,t,w,w-2)$-CPECC code, and let $G_i=(V_i,E_i)$ be a graph for every $i\in[M]$, where $V_i=\bigcup_{B\in\mathcal P_i}B$ and $E_i=\bigcup_{B\in\mathcal P_i}\{\{u,v\}:u,v\in B,u\neq v\}$. Since any two $w$-subsets from different $\mathcal{P}_i$'s have symmetric difference of size at least $2(w-2)+1$ by Property $\textbf{C}(w-2)$, any two $w$-subsets from different $\mathcal{P}_i$'s have at most one common point, and hence $E_1,E_2,\dots,E_M$ are pairwise disjoint. Therefore, $\sum_{i=1}^Me(G_i)\leq\binom{n}{2}.$ Establishing a lower bound on $e(G_i)$ is crucial for deriving upper bounds on $(n,t,w,w-2)$-CPECC codes. To derive tighter upper bounds on $(n,t,w,w-2)$-CPECC codes, Theorem \ref{thm:lb_e(G)_(n,t,w,w-2)} gives a lower bound on $e(G_i)$ and a stronger bound when $G_i$ is $K_{w+t}$-free. We first state the next two lemmas before proving Theorem \ref{thm:lb_e(G)_(n,t,w,w-2)}.

\begin{lemma}\label{lem:critical-cover}
Let $G$ be a graph and let $r\geq1$ be an integer. If $\tau_v(G)\geq r$, then $G$ contains a subgraph $J$ satisfying $\tau_v(J)=r$ and $|V(J)|\leq2r$.
\end{lemma}

\begin{proof}
We may choose, among all subgraphs of $G$ with vertex-cover number at least $r$, one with the minimum number of edges, and denote this subgraph by $J$. Since $\tau_v(G)\geq r$, there is such a subgraph $J$. As isolated vertices are irrelevant to the vertex-cover number, we may assume without loss of generality that $J$ has no isolated vertices.

For every $e\in E(J)$, the minimal choice of $J$ gives
\begin{align}\label{eq:J-minus-e-upper}
\tau_v(J-e)\leq r-1.
\end{align}
A vertex cover of $J-e$, together with either endpoint of $e$, is a vertex cover of $J$. Hence
\begin{align}\label{eq:one-edge-cover}
\tau_v(J)\leq\tau_v(J-e)+1.
\end{align}
Since $\tau_v(J)\geq r$, \eqref{eq:J-minus-e-upper} and \eqref{eq:one-edge-cover} force $\tau_v(J)=r$ and $\tau_v(J-e)=r-1$. Thus $J$ is $\tau_v$-critical. By Theorem \ref{thm:EG}, $|V(J)|\leq2r$.
\end{proof}

\begin{lemma}\label{lem:ordering}
Let $k\geq1$ be an integer, and let $G$ be a graph with $\delta(G)\geq k-1$. Then there exists a set $X\subseteq V(G)$ such that $G[X]$ is $k$-colorable and
$|X|\geq \sum_{v\in V(G)}\frac{k}{d_G(v)+1}$.
\end{lemma}

\begin{proof}
Let $V(G)$ be ordered uniformly at random, and let $X$ be the set of vertices with at most $k-1$ neighbors preceding them in this ordering.

Color the vertices of $G[X]$ one by one, following the inherited order. When a vertex is reached, at most $k-1$ of its already-colored neighbors lie in $X$. These neighbors forbid at most $k-1$ colors, leaving at least one of the $k$ colors available. Thus $G[X]$ is $k$-colorable.

Fix $v\in V(G)$. In the random ordering, each of the $d_G(v)+1$ vertices in $\{v\}\cup N_G(v)$ is equally likely to appear in any position relative to the others. The vertex $v$ has at most $k-1$ earlier neighbors exactly when it is among the first $k$ positions in this local ordering. Since $\delta(G)\geq k-1$, $d_G(v)+1\geq k$, and so
\begin{align*}
\Pr(v\in X)=\frac{k}{d_G(v)+1}.
\end{align*}
Let $\mathbf{1}_{\{v\in X\}}$ be the indicator random variable for the event $\{v\in X\}$. By linearity of expectation,
\begin{align*}
\mathbb E|X|=\sum_{v\in V(G)}\mathbb E\mathbf{1}_{\{v\in X\}}=\sum_{v\in V(G)}\frac{k}{d_G(v)+1}.
\end{align*}
Since there exists an ordering for which $|X|\geq\mathbb E|X|$, the desired conclusion is obtained.
\end{proof}

\begin{theorem}\label{thm:lb_e(G)_(n,t,w,w-2)}
Let $t\geq1$ and $w>t+3$ be integers, and let $G$ be a graph with $\delta(G)\geq w-1$. Suppose that
$G-S\text{ contains a copy of }K_w\ \text{for every }S\subseteq V(G)\text{ with }|S|=t$.
Then
\begin{align}\label{eq:ordinary-bound}
e(G)\geq\binom{w+t}{2},
\end{align}
and equality holds if and only if $G\cong K_{w+t}$.

Furthermore, if $G$ is $K_{w+t}$-free, then
\begin{align}\label{eq:strong-bound}
e(G)\geq B(w,t),
\end{align}
where $B(w,t)=\binom{w+\lfloor(t-1)/2\rfloor}{2}+\binom{w+\lceil(t-1)/2\rceil}{2}$. This stronger bound is sharp.
\end{theorem}

\begin{proof} To ensure smooth reading of the paper, we move the proof to Appendix \ref{app-A}.
\end{proof}

The following theorem shows an upper bound for $(n,t,w,w-2)$-CPECC codes.
		
\begin{theorem}\label{thm:ub_C'(n,1,w,w-2)}
Let $n,t,w$ be positive integers with $n\geq t+w$, $t\geq1$ and $w>t+3$. Then
$$C'(n,t,w,w-2)\leq \left\lfloor\frac{n}{w+t}\Big\lfloor\frac{n-1}{w+t-1}\Big\rfloor\right\rfloor+
\left\lfloor\frac{n(n-1)-\left\lfloor\frac{n}{w+t}\big\lfloor\frac{n-1}{w+t-1}\big\rfloor\right\rfloor(w+t)(w+t-1)}{2\left(\binom{w+\lfloor(t-1)/2\rfloor}{2}+\binom{w+\lceil(t-1)/2\rceil}{2}\right)}\right\rfloor.$$
\end{theorem}

\begin{proof}
Let $\mathcal B=\mathcal P_1\cup\mathcal P_2\cup\dots\cup\mathcal P_M$ be the set-theoretic representation of an $(n,t,w,w-2)$-CPECC code. Let $G_i=(V_i,E_i)$ be a graph for every $i\in[M]$, where $V_i=\bigcup_{B\in\mathcal P_i}B$ and $E_i=\bigcup_{B\in\mathcal P_i}\{\{u,v\}:u,v\in B,u\neq v\}$.

Without loss of generality, suppose that $G_i$ contains a copy of $K_{w+t}$ for all $i\in[b]$, while $G_i$ is $K_{w+t}$-free for all $b+1\leq i\leq M$. For each $i\in[b]$, choose a set $W_i\subseteq V_i$ such that $G_i[W_i]\cong K_{w+t}$. By Property $\textbf{C}(w-2)$, the edge sets $E_1,E_2,\dots,E_M$ are pairwise disjoint, and hence $G_i[W_i]$, $i\in[b]$, are pairwise edge-disjoint. That is, $([n],\{W_i:i\in[b]\})$ is an $(n,w+t,1)$-packing. It follows from Lemma \ref{lem:J_bound} that
$$b\leq \left\lfloor\frac{n}{w+t}\Big\lfloor\frac{n-1}{w+t-1}\Big\rfloor\right\rfloor.$$

Every vertex of $V_i$ belongs to some block in $\mathcal P_i$, and hence has degree at least $w-1$ in $G_i$. Moreover, by Property $\textbf{A}(t)$, $G_i-S$ contains a copy of $K_w$ for every $S\subseteq V_i$ with $|S|=t$. By Theorem \ref{thm:lb_e(G)_(n,t,w,w-2)}, $e(G_i)\geq\binom{w+\lfloor(t-1)/2\rfloor}{2}+\binom{w+\lceil(t-1)/2\rceil}{2}$ for all $b+1\leq i\leq M$, and so
$$b\binom{w+t}{2}+(M-b)\left(\binom{w+\lfloor(t-1)/2\rfloor}{2}+\binom{w+\lceil(t-1)/2\rceil}{2}\right)\leq\sum_{i=1}^M e(G_i)\leq\binom{n}{2}.$$
Thus
\begin{align*}
M&=b+(M-b)\leq b+\frac{{n\choose 2}-{w+t\choose 2}b}{\binom{w+\lfloor(t-1)/2\rfloor}{2}+\binom{w+\lceil(t-1)/2\rceil}{2}}\\
&=b\left(1-\frac{{w+t\choose 2}}{\binom{w+\lfloor(t-1)/2\rfloor}{2}+\binom{w+\lceil(t-1)/2\rceil}{2}}\right)+\frac{{n\choose 2}}{\binom{w+\lfloor(t-1)/2\rfloor}{2}+\binom{w+\lceil(t-1)/2\rceil}{2}}.
\end{align*}
Since $1-\frac{{w+t\choose 2}}{\binom{w+\lfloor(t-1)/2\rfloor}{2}+\binom{w+\lceil(t-1)/2\rceil}{2}}>0$ by (\ref{eq:B-strict}),
we have
\begin{align*}
M&\leq b\left(1-\frac{{w+t\choose 2}}{\binom{w+\lfloor(t-1)/2\rfloor}{2}+\binom{w+\lceil(t-1)/2\rceil}{2}}\right)+\frac{{n\choose 2}}{\binom{w+\lfloor(t-1)/2\rfloor}{2}+\binom{w+\lceil(t-1)/2\rceil}{2}}\\
&\leq \left\lfloor\frac{n}{w+t}\Big\lfloor\frac{n-1}{w+t-1}\Big\rfloor\right\rfloor+\frac{n(n-1)-\left\lfloor\frac{n}{w+t}\big\lfloor\frac{n-1}{w+t-1}\big\rfloor\right\rfloor(w+t)(w+t-1)}{2\left(\binom{w+\lfloor(t-1)/2\rfloor}{2}+\binom{w+\lceil(t-1)/2\rceil}{2}\right)}.
\end{align*}
Since $M$ is an integer, the desired conclusion is obtained.
\end{proof}

Remark that taking $t=1$ in Theorem \ref{thm:ub_C'(n,1,w,w-2)} yields \cite[Theorem 3.4]{ltz} for $w\geq5$. The following theorem provides a construction of optimal $(n,t,w,w-2)$-CPECC codes for $t\geq1$.

\begin{theorem}\label{cor-2}
Let $t\geq1$ and $w>t+3$ be integers with $t\equiv1\pmod2$. When $n$ is large enough, if $n(n-w-\frac{t-1}{2})\equiv0\pmod{(w+t)(w+t-1)}$, $n-w-\frac{t-1}{2}\equiv0\pmod{w+t-1}$, and $n\equiv0\pmod{w+\frac{t-1}{2}}$, then
$$C'(n,t,w,w-2)=\frac{n(n-w-\frac{t-1}{2})}{(w+t)(w+t-1)}+\Big\lfloor\frac{n}{2w+t-1}\Big\rfloor.$$
\end{theorem}

\begin{proof}
Set $g=w+\frac{t-1}{2}$, $k=w+t$ and $u=\frac{n}{g}$. It follows from Theorem \ref{thm:asy_dgg} that there exists a $k$-GDD of type $g^u$. Let $(X,\mathcal G, \mathcal B)$ be a $k$-GDD of type $g^u$, where $\mathcal G=\{G_1,G_2,\ldots,G_u\}$. Set
$$\mathcal C=\left\{{B\choose w}: B\in \mathcal B\right\}\cup \left\{{G_{2i-1}\choose w}\cup {G_{2i}\choose w}: 1\leq i\leq \lfloor u/2\rfloor\right\}.$$

Now we show that $\mathcal C$ is an $(n,t,w,w-2)$-CPECC code. Every set belonging to a member of $\mathcal C$ has size $w$, so Property $\textbf{B}(w)$ holds.

Let $S\subseteq X$ be any set of size $t$. For every block $B\in\mathcal B$, $|B\setminus S|\ge k-t\ge w$, and so there is an $A\in\binom Bw$ such that $A\cap S=\emptyset$. For any \(G_{2i-1}\) and \(G_{2i}\), $|S\cap G_{2i-1}|+|S\cap G_{2i}|\le t$. Thus at least one of two groups, say \(G_j\), satisfies $|S\cap G_j|\le\left\lfloor\frac t2\right\rfloor=\frac{t-1}{2}$. Consequently, $|G_j\setminus S|\ge g-\frac{t-1}{2}=w$. Hence there is an $A\in\binom{G_j}{w}$ such that $A\cap S=\emptyset$. Thus Property $\textbf{A}(t)$ holds.

The definition of $k$-GDD implies that distinct blocks intersect in at most one point, and every block meets each group in at most one point. Let $A$ and $A'$ be two $w$-subsets belonging to distinct parts. If $A\subseteq B\in \mathcal B$ and $A'\subseteq B'\in\mathcal B$, where $B\neq B'$, then $|A\cap A'|\le |B\cap B'|\le1$.
If $A\subseteq B\in \mathcal B$ and $A'\subseteq G_j\in\mathcal G$, then $|A\cap A'|\le |B\cap G_j|\le1$. If $A\subseteq G_i\in \mathcal G$ and $A'\subseteq G_j\in\mathcal G$, where $i\ne j$, then $A\cap A'=\varnothing$. To sum up, $|A\Delta A'|=|A|+|A'|-2|A\cap A'|\geq2w-3$. Thus Property $\textbf{C}(w-2)$ holds.

Next, we check the size of the $(n,t,w,w-2)$-CPECC code $\mathcal C$. The number of blocks of this GDD is
$$\frac{g^2u(u-1)}{k(k-1)}=\frac{n(n-g)}{(w+t)(w+t-1)}=\frac{n(n-w-\frac{t-1}{2})}{(w+t)(w+t-1)}.$$
Moreover, pairing the groups of the GDD yields $\lfloor u/2\rfloor=\lfloor n/(2w+t-1)\rfloor$ additional parts. Then
$$|\mathcal C|=\frac{n(n-w-\frac{t-1}{2})}{(w+t)(w+t-1)}+\Big\lfloor\frac{n}{2w+t-1}\Big\rfloor.$$
Since $1<w+\frac{t-1}{2}\leq w+t-1$, $n(n-w-\frac{t-1}{2})\equiv0\pmod{(w+t)(w+t-1)}$ and $n-w-\frac{t-1}{2}\equiv0\pmod{w+t-1}$, we have
$$\left\lfloor\frac{n}{w+t}\Big\lfloor\frac{n-1}{w+t-1}\Big\rfloor\right\rfloor=\left\lfloor\frac{n}{w+t}\cdot\frac{n-w-\frac{t-1}{2}}{w+t-1}
\right\rfloor=\frac{n(n-w-\frac{t-1}{2})}{(w+t)(w+t-1)},$$
and so the upper bound in Theorem \ref{thm:ub_C'(n,1,w,w-2)} is
$$\frac{n(n-w-\frac{t-1}{2})}{(w+t)(w+t-1)}+\left\lfloor\frac{n(n-1)-n(n-w-\frac{t-1}{2})}{4\binom{w+(t-1)/2}{2}}\right\rfloor
=\frac{n(n-w-\frac{t-1}{2})}{(w+t)(w+t-1)}+\Big\lfloor\frac{n}{2w+t-1}\Big\rfloor.$$
Thus $\mathcal C$ is an optimal $(n,t,w,w-2)$-CPECC code.
\end{proof}

\begin{theorem}\label{cor-gdd-1}
Let $w>6$ be an integer. If $n$ is large enough, $n(n-w-1)\equiv0\pmod{(w+2)(w+3)}$, $n-w-1\equiv0\pmod{w+2}$, and $n\equiv0\pmod{w+1}$, then
$$C'(n,3,w,w-2)=\frac{n(n-w-1)}{(w+3)(w+2)}+\Big\lfloor\frac{n}{2(w+1)}\Big\rfloor.$$
\end{theorem}

\begin{proof}
Apply Theorem \ref{cor-2} with $t=3$.
\end{proof}

\section{Optimal $(n,t,w,w-2)$-CPECC codes for $t\geq w-3$}

In \cite{ltz}, the authors established an upper bound for $(n,t,w,w-2)$-CPECC codes. However, this upper bound cannot always be attained. This section improves the upper bound in \cite{ltz} (see Theorem \ref{thm:ubj_(n,t,w,w-2)-CPECC} (1)) and constructs some new optimal families of $(n,t,w,w-2)$-CPECC codes.

\begin{theorem}{\label{thm:ubj_(n,t,w,w-2)-CPECC}}
Let $n,t$ and $w$ be positive integers with $n\geq t+w$.
\begin{enumerate}
   \item[$(1)$] {\rm \cite{ltz}} If $w\geq3$, then $C'(n,t,w,w-2)\leq\frac{{n\choose 2}}{(2w-3)(t+1)}.$
   \item[$(2)$] {\rm \cite{z-z}} If $w\geq t^2+2t+2\geq 5$, then $C'(n,t,w,w-2)\leq \frac{\binom{n}{2}}{\binom{w+t}{2}}.$
\end{enumerate}
\end{theorem}

\subsection{Upper bound for $(n,t,w,w-2)$-CPECC codes}

We employ the same idea as in Section \ref{sec:t<w-3} to establish an upper bound for $(n,t,w,w-2)$-CPECC codes. To obtain an upper bound for $(n,t,w,w-2)$-CPECC codes better than that in \cite{ltz}, we first introduce some definitions and results.

A {\em hypergraph} is a pair $H=(V(H),\mathcal E(H))$, where $V(H)$ is a finite vertex set and $\mathcal E(H)\subseteq2^{V(H)}$ is a family of subsets called {\em hyperedges}. It is {\em $r$-uniform} if every hyperedge has exactly $r$ vertices. A set $T\subseteq V(H)$ is a {\em transversal} of $H$ if $T\cap E\neq\varnothing$ for every $E\in\mathcal E(H)$. The {\em transversal number} is $\tau(H)=\min\{|T|:T\subseteq V(H)$ is a transversal of $H\}$. Note that $\tau(H)=0$ for an empty hypergraph $H$. A set $S\subseteq V(H)$ is {\em independent} in $H$ if no hyperedge of $H$ is contained in $S$. Note that if $S$ is independent in $H$, then $V(H)\setminus S$ is a transversal of $H$. For a hypergraph $H$, a {\em $2$-shadow} of $H$ is a graph $Q=(V(Q),E(Q))$, where $V(Q)=\bigcup_{E\in \mathcal E(H)}E$ and $E(Q)=\{\{x,y\}:\text{there is an }E\in\mathcal E(H) \text{ such that }x,y\in E,\ x\neq y\}$. Note that if $H$ is $(d+1)$-uniform, each hyperedge of $H$ induces a copy of $K_{d+1}$ in $Q$. For $v\in V(H)$, define $H-v=(V(H)\setminus\{v\},\{E\in\mathcal E(H):v\notin E\})$.

\begin{lemma}\label{lem:random}
Let $d\geq1$ be an integer, let $H$ be a $(d+1)$-uniform hypergraph with $V(H)=[n]$, and let $Q$ be the $2$-shadow of $H$. Suppose that $m=\tau(H)$ and $b=2d-1$. Then
\begin{enumerate}
\item[$(a)$] $m\leq\sum_{v\in V(Q)}\frac{d_Q(v)-d+1}{d_Q(v)+1}$, and
\item[$(b)$] $e(Q)-bm\geq\sum_{v\in V(Q)}\frac{(d_Q(v)-(2d-2))(d_Q(v)-(2d-1))}{2(d_Q(v)+1)}\geq0$.
\end{enumerate}
\end{lemma}

\begin{proof}
$(a)$ Choose a uniformly random permutation $\pi$ of $[n]$. For a vertex $v$, let $\pi(v)$ denote its position in the permutation. Define
$$ S_\pi=\left\{v\in V(Q):\bigl|\{u\in N_Q(v):\pi(u)<\pi(v)\}\bigr|\leq d-1\right\}. $$

We first show that the induced graph $Q[S_\pi]$ is $(d-1)$-degenerate. For every nonempty set $U\subseteq S_\pi$, choose $v\in U$ whose position $\pi(v)$ is largest among the vertices of $U$. Then every vertex of $U\setminus\{v\}$ appears earlier than $v$ in the permutation. Since $v\in S_\pi$, $v$ has at most $d-1$ neighbours among all vertices earlier than it, and so $d_{Q[U]}(v)\leq d-1.$
Hence the graph $Q[S_\pi]$ is $(d-1)$-degenerate. By \cref{lem:degenerate-colour}, $Q[S_\pi]$ is $d$-colourable.

Assume that $S_\pi$ contains some hyperedge $E\in\mathcal E(H)$. By the definition of the $2$-shadow, the graph $Q[E]$ is a copy of $K_{d+1}$. But $K_{d+1}$ is not $d$-colourable, a contradiction with the $d$-colourability of $Q[S_\pi]$. Hence $S_\pi$ contains no hyperedge of $H$ and is independent in $H$.
By the definitions of an independent set and a transversal, $V(Q)\setminus S_\pi$ is a transversal, and so $|V(Q)|-|S_\pi|\geq m$. Taking expectations yields $$|V(Q)|-m\geq\E|S_\pi|.$$

We now calculate the expectation of $|S_\pi|$. Fix a vertex $v$. Among the $d_Q(v)+1$ vertices in $N_Q[v]=N_Q(v)\cup\{v\}$, the relative position of $v$ is uniformly distributed over $1,2,\dots,d_Q(v)+1$. The event $v\in S_\pi$ occurs exactly when at most $d-1$ neighbours precede $v$, equivalently, when $v$ is among the first $d$ vertices of $N_Q[v]$ in the random order. Thus
$$ \Prb(v\in S_\pi)=\frac{d}{d_Q(v)+1}. $$
By linearity of expectation,
$$ \E|S_\pi|=\sum_{v\in V(Q)}\frac{d}{d_Q(v)+1}. $$
Therefore,
\begin{align*}
m \leq |V(Q)|-\sum_{v\in V(Q)}\frac{d}{d_Q(v)+1} =\sum_{v\in V(Q)}\left(1-\frac{d}{d_Q(v)+1}\right) =\sum_{v\in V(Q)}\frac{d_Q(v)-d+1}{d_Q(v)+1}.
\end{align*}
$($b$)$ Let $ X=\sum_{v\in V(Q)}\frac{d_Q(v)-d+1}{d_Q(v)+1}$ and $D=e(Q)-bm$. By $(a)$, $D=e(Q)-bX+b(X-m)\geq e(Q)-bX$.
Note that
\begin{align}
D\geq e(Q)-bX=&\sum_{v\in V(Q)}\left(\frac{d_Q(v)}{2}-\frac{b(d_Q(v)-d+1)}{d_Q(v)+1}\right)\notag\\
&=\sum_{v\in V(Q)}\frac{(d_Q(v)-(2d-2))(d_Q(v)-(2d-1))}{2(d_Q(v)+1)}\geq0,\notag
\end{align}
where the first equality holds since $2e(Q)=\sum_v d_Q(v)$ and the last inequality holds because $d_Q(v)$ is an integer and $2d-2$ and $2d-1$ are consecutive integers.
\end{proof}

Let $\Z$ be the set of integers and $\Z_{\geq0}=\{x\in\Z:x\geq0\}$. Fix an integer $d\geq3$. For each $p\in\Z_{\geq0}$, define the integer interval
$$I_p=[p(d-1),pd]_{\Z}=\{p(d-1),p(d-1)+1,\dots,pd\}.$$
In particular, $I_0=\{0\}$. For an integer $x$ and a nonempty set $S\subseteq\Z$, define $\dist(x,S)=\min_{y\in S}|x-y|$. For $m\in\Z_{\geq0}$, define
$$z_d(m)=\min_{p\in\Z_{\geq0}}\dist(m,I_p)=\dist\left(m,\bigcup_{p\geq0}I_p\right).$$

\begin{lemma}\label{lem:z-properties}
For all $x,y\in\Z_{\geq0}$, $z_d(x+y)\leq z_d(x)+z_d(y)$ and $|z_d(x+1)-z_d(x)|\leq1.$
\end{lemma}

\begin{proof}
Let $p,q\in\Z_{\geq0}$ and $u\in I_p$, $v\in I_q$ such that $|x-u|=z_d(x)$ and $|y-v|=z_d(y)$. Let $I_p+I_q=\{r+s:r\in I_p,\ s\in I_q\}$. Since
$I_p+I_q=I_{p+q}$, $u+v\in I_{p+q}$.
Therefore
\begin{align*}
z_d(x+y)\leq|x+y-(u+v)|\leq|x-u|+|y-v|=z_d(x)+z_d(y).
\end{align*}

Set $S=\bigcup_{p\geq0}I_p$. Choose $u\in S$ with $|x-u|=z_d(x)$. Then
$$z_d(x+1)\leq|x+1-u|\leq|x-u|+1=z_d(x)+1.$$
Conversely, choose $v\in S$ with $|x+1-v|=z_d(x+1)$. Then
$$z_d(x)\leq|x-v|\leq|x+1-v|+1=z_d(x+1)+1.$$
Combining the above two inequalities yields $|z_d(x+1)-z_d(x)|\leq1$.
\end{proof}

The following theorem shows a sharper estimate $D\geq z_d(m)$, where $b=2d-1$ and $D=e(Q)-bm$.

\begin{theorem}\label{thm:shadow}
Let $d\geq3$ and $b=2d-1$. If $H$ is a $(d+1)$-uniform hypergraph with $m=\tau(H)$ and $Q$ is the $2$-shadow of $H$, then
$$e(Q)\geq bm+z_d(m).$$
\end{theorem}

\begin{proof} To ensure smooth reading of the paper, we move the proof to Appendix \ref{app-B}.
\end{proof}

Given a graph $G$, define a $w$-uniform hypergraph $\mathcal K_w(G)=\bigl(V(G),\mathcal E_w(G)\bigr),$ where
$$\mathcal E_w(G)=\{E\subseteq V(G):|E|=w\text{ and }G[E]\cong K_w\}.$$
Note that the hyperedges of $\mathcal K_w(G)$ are exactly the vertex sets of copies of $K_w$ in $G$.

\begin{theorem}\label{lem-4.1}
Let $w\geq4$ and $t\geq w-3$ be integers, and let $G$ be a graph with $|V(G)|\geq t$ and $\delta(G)\geq w-1$. If $G-S$ contains a copy of $K_w$ for every $S\subseteq V(G)$ with $|S|=t$, then
$$ e(G)\geq (2w-3)(t+1)+z_{w-1}(t+1).$$
\end{theorem}

\begin{proof}
We claim that $\tau(\mathcal K_w(G))\geq t+1$. Indeed, suppose, to the contrary, that $\mathcal K_w(G)$ has a transversal $T\subseteq V(G)$ with $|T|\leq t$. Because $|V(G)|\geq t$, we enlarge $T$ to a set $T^*\subseteq V(G)$ with $T\subseteq T^*$ and $|T^*|=t$. Since $T$ meets the vertex set of every copy of $K_w$ in $G$, the larger set $T^*$ also meets every copy of $K_w$. Consequently, $G-T^*$ contains no $K_w$, contradicting the hypothesis. 

Because $G$ is finite, the hypergraph $\mathcal K_w(G)$ has finitely many hyperedges. Enumerate them as $E_1,E_2,\dots,E_N$. Let $H_0$ be the empty hypergraph, and for $1\leq i\leq N$, define
$$H_i=\bigl(V(G),\{E_1,\dots,E_i\}\bigr).$$
Therefore, $\tau(H_0)=0$ and $\tau(H_N)=\tau(\mathcal K_w(G))\geq t+1$ by the above claim.

For $1\leq i\leq N$, let $T$ be a minimum transversal of $H_{i-1}$. Choose any vertex $x\in E_i$. Then $T\cup\{x\}$ meets every hyperedge of $H_i$, so $\tau(H_i)\leq\tau(H_{i-1})+1$. Let $i$ be the smallest index for which $\tau(H_i)\geq t+1$. By the choice of $i$, $\tau(H_{i-1})\leq t$, and so $\tau(H_i)\leq\tau(H_{i-1})+1\leq t+1$. Thus $\tau(H_i)=t+1$.

Let $H=H_i$ and $Q$ be a 2-shadow of $H$. Note that  $e(G)\geq e(Q)$. Applying \cref{thm:shadow} with $d=w-1$, $e(Q)\geq (2w-3)(t+1)+z_{w-1}(t+1)$. Thus $e(G)\geq (2w-3)(t+1)+z_{w-1}(t+1)$.
\end{proof}

\begin{theorem}[Upper Bound]{\label{thm:upb_(n,t,w,w-2)-CPECC}}
Let $n,t,w$ be positive integers with $n\geq t+w$, $w\geq4$, and $t\geq w-3$. Then
$$C'(n,t,w,w-2)\leq\frac{n(n-1)}{2((2w-3)(t+1)+z_{w-1}(t+1))}.$$
\end{theorem}

\begin{proof}
Let $\mathcal{C}={\cal C}_1\cup{\cal C}_2\cup\ldots\cup{\cal C}_M$ be an $(n,t,w,w-2)$-CPECC code and ${\cal B}={\cal P}_1\cup\ldots\cup{\cal P}_M$ the corresponding set-theoretic representation of $\mathcal{C}$ over $[n]$. Let $G_i=(V_i,E_i)$ be a graph for every $i\in[M]$, where $V_i=\bigcup_{B\in\mathcal P_i}B$ and $E_i=\bigcup_{B\in\mathcal P_i}\{\{u,v\}:u,v\in B,u\neq v\}$.

Since any two codewords from different $\mathcal{C}_i$'s have Hamming distance at least $2(w-2)+1$, any two $w$-subsets from different $\mathcal{P}_i$'s have at most one common point, and hence $E_1,E_2,\dots,E_M$ are pairwise disjoint. Thus
$$M\leq\frac{\binom{n}{2}}{\min_{i\in[M]}e(G_i)}.$$

Every vertex of $V_i$ belongs to some block in $\mathcal P_i$, and hence $\delta(G_i)\geq w-1$. Moreover, Property $\textbf{A}(t)$ implies that $|V_i|\geq t+w$. Indeed, if $|V_i|\leq t$, then a $t$-subset of $[n]$ can be chosen to contain $V_i$, contradicting Property $\textbf{A}(t)$. If $t<|V_i|<t+w$, choose a $t$-subset $S\subseteq V_i$; then $|V_i\setminus S|<w$, again contradicting Property $\textbf{A}(t)$. Furthermore, by Property $\textbf{A}(t)$, $G_i-S$ contains a copy of $K_w$ for every $S\subseteq V_i$ with $|S|=t$. It follows from Theorem \ref{lem-4.1} that $e(G_i)\geq (2w-3)(t+1)+z_{w-1}(t+1)$ for each $i\in[M]$. Thus, $$M\leq\frac{\binom{n}{2}}{(2w-3)(t+1)+z_{w-1}(t+1)}=\frac{n(n-1)}{2((2w-3)(t+1)+z_{w-1}(t+1))}.$$
This completes the proof.
\end{proof}

\subsection{Restatement of $(2w-3)(t+1)+z_{w-1}(t+1)$}
Since computing $z_{w-1}(t+1)$ is difficult, this subsection presents a new expression for $(2w-3)(t+1)+z_{w-1}(t+1)$.

Let $\Z_{>0}=\{x\in \Z:x>0\}$ and $A=\{a_1,\dots,a_k\}\subseteq\Z_{>0}$ with $\gcd(a_1,\dots,a_k)=1$. Its nonnegative integer semigroup is
$\langle A\rangle=\left\{\sum_{i=1}^k n_i a_i:n_i\in\Z_{\geq0}\right\}.$ The \emph{Frobenius number} of $A$ is
$$ F(A) = \left\{
\begin{array}{cl}
\max\bigl(\Z_{\geq0}\setminus\langle A\rangle\bigr) & \text{if } 1 \notin A,\\
-1 & \text{if } 1\in A.
\end{array} \right. $$
Let $J_p=[pa,p(a+3)]_{\Z}$ for $p\in\Z_{\geq0}$.

\begin{lemma}\label{lem:Jp}
For any integer $a\geq1$, using exactly $p$ elements from $\{a,a+1,a+2,a+3\}$ (with repetition allowed), the set of all possible sums is precisely $J_p$, where $J_p=[pa,p(a+3)]_{\Z}$ for $p\in\Z_{\geq0}$.
\end{lemma}

\begin{proof}
Write each chosen denomination uniquely as $a+\varepsilon_i$ with $\varepsilon_i\in\{0,1,2,3\}$. A sum using exactly $p$ coins has the form $pa+\sum_{i=1}^{p}\varepsilon_i$, where $0\leq\sum_{i=1}^{p}\varepsilon_i\leq3p$. Thus every such sum belongs to $J_p$.

Let $u\in\{0,1,\dots,3p\}$. Write $u$ by Euclidean division by $3$: $u=3q+r$, where $q\in\Z_{\geq0}$ and $0\leq r\leq2$. If $r=0$, use $q$ entries equal to $3$ and all remaining entries equal to $0$. If $r>0$, use $q$ entries equal to $3$, one entry equal to $r$, and all remaining entries equal to $0$. Because $u\leq3p$, the number of nonzero entries needed is at most $p$. Therefore every $u\in[0,3p]_{\Z}$ is a sum of $p$ elements of $\{0,1,2,3\}$, and every integer in $J_p$ is representable using exactly $p$ coins.
\end{proof}

\begin{lemma}\label{prop:F}
For any integer $a\geq1$, $F(\{a,a+1,a+2,a+3\})=a\left\lfloor\frac{a+1}{3}\right\rfloor-1.$
\end{lemma}

\begin{proof}
By Lemma \ref{lem:Jp}, sums using exactly $p$ coins form the interval $J_p$. The consecutive intervals $J_p$ and $J_{p+1}$ have no integer gap between them exactly when $(p+1)a\leq p(a+3)+1,$ which is equivalent to $3p+1\geq a.$
The smallest nonnegative integer satisfying this inequality is $$p_0=\left\lfloor\frac{a+1}{3}\right\rfloor.$$
Thus $J_{p_0},J_{p_0+1},\dots$ join without gaps and cover every integer from $ap_0$ onward. Hence every integer $x\geq ap_0$ is representable.

It remains to show that $ap_0-1$ is not representable. For $a=1$, $p_0=0$. Since $\bigcup_{j\geq0}J_j=\Z_{\geq0}$, $ap_0-1=-1$ is the largest integer that is not representable. Suppose that $a\geq2$. If $p\leq p_0-1$, then $$p(a+3)\leq(p_0-1)(a+3)=ap_0+(3p_0-a-3)\leq ap_0-2,$$ where the last inequality holds by $3p_0\leq a+1$.
Thus no interval $J_p$ with $p\leq p_0-1$ reaches $ap_0-1$. On the other hand, if $p\geq p_0$, then every element of $J_p$ is at least $pa\geq ap_0$.
Therefore $ap_0-1$ is not representable and is the largest nonrepresentable nonnegative integer.
\end{proof}

Call $a$ and $a+3$ the \emph{endpoint denominations}, and call $a+1$ and $a+2$ the \emph{middle denominations}. For every representable $x\in\Z_{\geq0}$, define
$$c(x)=\min\{l_1+l_4:x=l_1a+l_2(a+1)+l_3(a+2)+l_4(a+3),\ l_1,l_2,l_3,l_4\in\Z_{\geq0}\}.$$
Note that $c(x)$ is the smallest possible number of endpoint coins in a representation of $x$.

\begin{lemma}\label{lem:lambda_p(x)}
$c(x)=\min_{\substack{p\geq0, x\in J_p}}\dist(x,I'_p)$, where $I'_p=[p(a+1),p(a+2)]_{\Z}$.
\end{lemma}

\begin{proof}
Fix $p\geq0$ and $x\in J_p$. Write uniquely $x=pa+u$, where $0\leq u\leq3p$. Let $\lambda_p(x)$ be the minimum number of endpoint coins among all representations of $x$ using exactly $p$ coins.

We claim that $\lambda_p(x)=\dist\bigl(x,[p(a+1),p(a+2)]_{\Z}\bigr)=\dist(x,I'_p)$. Indeed, represent each coin as $a+\varepsilon_i$ with $\varepsilon_i\in\{0,1,2,3\}$. Endpoint coins correspond exactly to $\varepsilon_i\in\{0,3\}$.

(1) For $0\leq u\leq p$, let $k$ be the number of endpoint coins in a representation using $p$ coins. The remaining $p-k$ coins are middle coins, so each of their $\varepsilon_i$ values is at least $1$. Therefore, $u=\sum_i\varepsilon_i\geq p-k$, which gives $k\geq p-u$. Equality is attained by taking $p-u$ values equal to $0$ and $u$ values equal to $1$. Hence $\lambda_p(x)=p-u$.

(2) For $p\leq u\leq2p$, take $u-p$ values equal to $2$ and $2p-u$ values equal to $1$. Their number is $(u-p)+(2p-u)=p$, and their sum is $2(u-p)+(2p-u)=u$. No endpoint value $0$ or $3$ is used, so $\lambda_p(x)=0$.

(3) Suppose that $2p\leq u\leq3p$. Let $k$ be the number of endpoint coins in a representation using $p$ coins. Note that every $\varepsilon_i$ of endpoint coins is at most $3$. The remaining $p-k$ coins are middle coins, so each of their $\varepsilon_i$ values is at most $2$. Then $u\leq3k+2(p-k)=2p+k$. Thus $k\geq u-2p$. Equality is attained by taking $u-2p$ values equal to $3$ and $3p-u$ values equal to $2$. Hence $\lambda_p(x)=u-2p$.

To sum up, the above three formulas for $\lambda_p(x)$ give exactly the distance from $x=pa+u$ to the interval $[p(a+1),p(a+2)]_{\Z}$.

By the definition of $c(x)$, $c(x)=\min_{\substack{p\geq0, x\in J_p}}\dist(x,I'_p)$.
\end{proof}

\begin{lemma}\label{lem:projection}
For $w\geq4$, let $d=w-1$ and $a=w-3$. If $x\geq a\lfloor\frac{a+1}{3}\rfloor$, then $z_d(x)=c(x)$.
\end{lemma}

\begin{proof}
Since $a+1=d-1$, $I_p=I'_p$ for each $p\in\Z_{\geq0}$. Lemma \ref{lem:lambda_p(x)} implies $c(x)=\min_{\substack{p\geq0, x\in J_p}}\dist(x,I_p)$. Recall $z_d(x)=\min_{p\in\Z_{\geq0}}\dist(x,I_p)$. To prove $z_d(x)=c(x)$, it suffices to prove $c(x)\leq \dist(x,I_p)$  for every $p\geq0$ with $x\not\in J_p$.

We first show that $c(x)\leq R$ for each $x\geq a\left\lfloor\frac{a+1}{3}\right\rfloor$, where $R=\left\lceil\frac a3\right\rceil$. For each $p\geq0$, define
$$B_p=J_p\cap[p(a+1)-R,\,p(a+2)+R]_{\Z}.$$
We claim that $B_R,B_{R+1},\dots$ join without integer gaps. Indeed, for $p\geq R$, $p(a+1)-R\geq pa$ and $p(a+2)+R\leq p(a+3)$. Thus $B_p=[p(a+1)-R,\,p(a+2)+R]_{\Z}.$
Since the integer gap between $B_p$ and $B_{p+1}$ is $\bigl((p+1)(a+1)-R\bigr)-\bigl(p(a+2)+R\bigr)-1=a-p-2R\leq a-3R\leq0$, we know that $B_R,B_{R+1},\dots$ join without integer gaps.

Let $p_0=\left\lfloor\frac{a+1}{3}\right\rfloor$.
If $a\equiv1\pmod3$, then $p_0=R-1$. Since $p_0(a+1)-R=p_0a-1$ and $p_0(a+2)+R=p_0(a+3)+1$, $B_{p_0}=J_{p_0}$. Since $p_0=\frac{a-1}{3}$ and $R=\frac{a+2}{3}$,
$$\max B_{p_0}+1=p_0(a+3)+1=\frac{a^2+2a}{3}=Ra=\min B_R.$$
Combined with the above claim, $B_{p_0},B_{p_0+1},\dots$ join without integer gaps. If $a\equiv0$ or $2\pmod3$, then $p_0=R$, and so the above claim implies $B_{p_0},B_{p_0+1},\dots$ join without integer gaps. Since the left endpoint of $B_{p_0}$ is $ap_0$, every integer $x\geq ap_0$ belongs to some $B_p$. Note that $B_p$ consists exactly of the integers that are representable using $p$ coins and are at distance at most $R$ from $I_p$. Therefore for every $x\geq ap_0$, there exists a \emph{legal} $p$, meaning $x\in J_p$, such that $\dist(x,I_p)\leq R.$ By Lemma \ref{lem:lambda_p(x)}, we have
$$c(x)\leq R.$$

It remains to show that every \emph{illegal} value of $p$, meaning $x\notin J_p$, satisfies
\begin{equation}\label{eq:illegal-R}
\dist(x,I_p)\geq R.
\end{equation}
There are only two cases: $x>p(a+3)$ and $x<pa$. Writing $a$ as $3u$, $3u+1$, or $3u+2$ gives, respectively,
\begin{equation}\label{eq:a=3u+123}
\begin{array}{c|c|c|c}
a&p_0&R&a+2-2p_0\\ \hline
3u&u&u&u+2\\
3u+1&u&u+1&u+3\\
3u+2&u+1&u+1&u+2
\end{array}
\end{equation}

\emph{$(1)$ $x>p(a+3)$.}
Since $x$ is an integer, we have $x\geq p(a+3)+1$. As $x$ lies to the right of $I_p$, we have
$$\dist(x,I_p)=x-p(a+2)\geq p+1.$$
If $p\geq R$, then $\dist(x,I_p)\geq R$. If $p_0\leq p<R$, then $a\equiv1\pmod3$ and $p=p_0=R-1$. Hence $\dist(x,I_p)\geq p+1=R$. For $p\leq p_0-1$, since $x\geq ap_0$, we have
\begin{align*}
\dist(x,I_p)
&\geq ap_0-(p_0-1)(a+2)
=a+2-2p_0\geq R,
\end{align*}
where the last inequality holds by \eqref{eq:a=3u+123}.

\emph{$(2)$ $x<pa$.}
Since $x\geq ap_0$, $p\geq p_0+1$. Since $x\leq pa-1$, $x$ lies to the left of $I_p$, and so
$$\dist(x,I_p)=p(a+1)-x\geq p+1\geq p_0+2\geq R.$$
Thus \eqref{eq:illegal-R} holds for every illegal $p$. This completes the proof.
\end{proof}

Define
$$\mathcal A_{w,t}=\{(l_1,l_2,l_3,l_4)\in\mathbb Z_{\geq0}^4:l_1(w-3)+l_2(w-2)+l_3(w-1)+l_4w=t+1\},$$
and, whenever $\mathcal A_{w,t}\neq\varnothing$, define
\begin{equation}\label{def:e_{w,t}}
e_{w,t}=\min_{(l_1,l_2,l_3,l_4)\in\Awt}
\Biggl[
\binom{2w-4}{2}l_1+
\binom{2w-3}{2}l_2+
\binom{2w-2}{2}l_3+
\binom{2w-1}{2}l_4
\Biggr].
\end{equation}
The following theorem provides a new expression for $(2w-3)(t+1)+z_{w-1}(t+1)$.

\begin{theorem}\label{thm:e_wt=bs+zds}
For $w\geq4$, let $d=w-1$, $b=2d-1$. If $t\geq (w-3)\lfloor\frac{w-2}{3}\rfloor-1$, then
\begin{equation}\label{eq:firstcandidate}
\begin{aligned}
b(t+1)+z_d(t+1)=e_{w,t}.
\end{aligned}
\end{equation}
\end{theorem}

\begin{proof}
Since $t\geq (w-3)\lfloor\frac{w-2}{3}\rfloor-1$, we have $\mathcal A_{w,t}\not=\emptyset$ by \cref{prop:F}. Applying \cref{lem:projection} to $x=t+1$ and $a=w-3$ gives
\begin{equation}
z_d(t+1)=c(t+1)=\min_{(l_1,l_2,l_3,l_4)\in\mathcal A_{w,t}}(l_1+l_4).\notag
\end{equation}
Note that
\begin{equation}
\binom{d+q}{2}=\frac{(q-d+b)(q-d+1+b)}{2}=bq+\frac{(q-d)(q-d+1)}2.\notag
\end{equation}
Applying the above identity to $q=d-2,d-1,d$ and $d+1$ gives
\begin{equation}
\begin{aligned}
\binom{2w-4}{2}&=b(w-3)+1,\ \ \
\binom{2w-3}{2}=b(w-2),\\
\binom{2w-2}{2}&=b(w-1),\text{ and }\
\binom{2w-1}{2}=bw+1.\notag
\end{aligned}
\end{equation}
For $(l_1,l_2,l_3,l_4)\in\Awt$,
$$
\begin{aligned}
&\binom{2w-4}{2}l_1+\binom{2w-3}{2}l_2+\binom{2w-2}{2}l_3+\binom{2w-1}{2}l_4\\
&=b\bigl(l_1(w-3)+l_2(w-2)+l_3(w-1)+l_4w\bigr)+l_1+l_4\\
&=b(t+1)+l_1+l_4.
\end{aligned}
$$
Taking the minimum over all $(l_1,l_2,l_3,l_4)\in\Awt$, the desired conclusion is obtained.
\end{proof}

By Theorems \ref{thm:upb_(n,t,w,w-2)-CPECC} and \ref{thm:e_wt=bs+zds}, we have the following theorem which shows a new expression of an upper bound for $(n,t,w,w-2)$-CPECC codes.
\begin{theorem}{\label{thm:upb_(n,t,w,w-2)-CPECC_new}}
Let $n,t,w$ be positive integers with $n\geq t+w$, $w\geq4$, and $t\geq \max\{w-3,(w-3)\lfloor\frac{w-2}{3}\rfloor-1\}$. Then
$$C'(n,t,w,w-2)\leq\frac{n(n-1)}{2e_{w,t}},$$
where $e_{w,t}$ is defined in $(\ref{def:e_{w,t}})$.
\end{theorem}

\begin{remark}\label{remark-1}
$(i)$ Suppose that $w\geq4$, $t\geq\max\{w-3,(w-3)\lfloor\frac{w-2}{3}\rfloor-1\}$, and that there exist nonnegative integers $l_2,l_3$ satisfying $t+1=l_2(w-2)+l_3(w-1)$. Then $${2w-3\choose2}l_2+{2w-2\choose2}l_3=(2w-3)(t+1),$$
and so $e_{w,t}\leq (2w-3)(t+1)$. By \eqref{eq:firstcandidate}, $e_{w,t}\geq b(t+1)=(2w-3)(t+1).$ Hence $e_{w,t}=(2w-3)(t+1)$. Theorems \ref{thm:upb_(n,t,w,w-2)-CPECC} and \ref{thm:e_wt=bs+zds} imply
$$C'(n,t,w,w-2)\leq\frac{n(n-1)}{2(2w-3)(t+1)},$$
which coincides with the bound given in \cite{ltz}.

$(ii)$ Let $w\geq4$ and $t\geq\max\{(w-3)\lfloor\frac{w-2}{3}\rfloor-1,w-3\}$. Suppose that $t+1=(w-3)l$ for some positive integer $l$ satisfying $3l\leq w-1$. Consider any representation of $t+1$ by elements of $\{w-3,w-2,w-1,w\},$ and suppose that it uses $p$ terms. Since every term is at least $w-3$ and $t+1=(w-3)l$, $p(w-3)\leq t+1=(w-3)l$, and so $$p\leq l.$$ On the other hand, $(l-1)w<(w-3)l=t+1$ by $3l\leq w-1$. Therefore $t+1$ cannot be represented using at most $l-1$ terms, that is, $$p\geq l.$$
Thus every representation uses exactly $l$ terms. Since each of these $l$ terms is at least $w-3$ and their sum is exactly $l(w-3)$, all of them must be equal to $w-3$. Consequently, $c(t+1)=l$.
By Lemma \ref{lem:projection}, $z_{w-1}(t+1)=l$. Therefore, by Theorem~\ref{thm:upb_(n,t,w,w-2)-CPECC}, we have
$$C'(n,t,w,w-2)\leq\frac{n(n-1)}{2((2w-3)(t+1)+l)}.$$
Note that the above upper bound is strictly smaller than the bound given in \cite{ltz}.
\end{remark}

\subsection{A construction of optimal $(n,t,w,w-2)$-CPECC codes}
In this section, we construct a new family of optimal $(n,t,w,w-2)$-CPECC codes whose size attains the upper bound in $(ii)$ of Remark \ref{remark-1}.
\begin{theorem}\label{new-22}
Let $w\geq4$ and $t\geq\max\{(w-3)\lfloor\frac{w-2}{3}\rfloor-1,w-3\}$. Suppose that $t+1=(w-3)l$ for some positive integer $l$ satisfying $3l\leq w-1$ and $l(2w-4)\mid n$. If there exists an $(n,2w-4,1)$-RBIBD, then
$$C'(n,t,w,w-2)=\frac{n(n-1)}{2((2w-3)(t+1)+l)}.$$
\end{theorem}

\begin{proof}
By Remark \ref{remark-1} $(ii)$, $$C'(n,t,w,w-2)\leq\frac{n(n-1)}{2((2w-3)(t+1)+l)}.$$

Let $(X,\mathcal B)$ be an $(n,2w-4,1)$-RBIBD. Its block set admits a resolution into parallel classes
$$\mathcal P_i=\left\{B_{i,h}:1\leq h\leq\frac{n}{2w-4}\right\},\ \text{ where }\ 1\leq i\leq\frac{n-1}{2w-5}.$$
Since $l(2w-4)\mid n$, each parallel class can be partitioned into groups of $l$ blocks. For a set $B$, write $\binom{B}{w}=\{W\subseteq B:|W|=w\}$. For any $1\leq i\leq\frac{n-1}{2w-5}$ and $1\leq j\leq\frac{n}{l(2w-4)}$, define
$$\mathcal P_{i,j}=\bigcup_{r=0}^{l-1}\binom{B_{i,jl-r}}{w}.$$
Let
$$\mathcal C=\bigcup_{i=1}^{\frac{n-1}{2w-5}}\bigcup_{j=1}^{\frac{n}{l(2w-4)}}\mathcal P_{i,j}.$$
Note that the size of an $(n,t,w,w-2)$-CPECC code $\mathcal C$ is $$\frac{n-1}{2w-5}\cdot\frac{n}{l(2w-4)}=\frac{n(n-1)}{2l(w-2)(2w-5)}=\frac{n(n-1)}{2((2w-3)(t+1)+l)}.$$

It remains to show that $\mathcal C$ is an $(n,t,w,w-2)$-CPECC code. Every member of $\mathcal C$ has size $w$, so Property $\textbf{B}(w)$ holds.

Let $T\subseteq X$ with $|T|=t$, and fix a part $\mathcal P_{i,j}$. The $l$ blocks $B_{i,jl},B_{i,jl-1},\ldots,B_{i,jl-l+1}$ belong to the same parallel class and hence are pairwise disjoint. Suppose that each of them contains at most $w-1$ points outside $T$. Since each block has size $2w-4$, every one of these $l$ blocks then meets $T$ in at least $(2w-4)-(w-1)=w-3$ points. Their pairwise disjointness gives $|T|\geq l(w-3)=t+1$, a contradiction. Hence at least one of these blocks contains at least $w$ points outside $T$, and therefore $\mathcal P_{i,j}$ contains a $w$-subset disjoint from $T$. Thus Property $\textbf{A}(t)$ holds.

Consider two $w$-subsets belonging to distinct parts. If their underlying RBIBD blocks lie in the same parallel class, those blocks are disjoint. If they lie in different parallel classes, then they intersect in at most one point by the definition of RBIBDs. Hence any two $w$-subsets from different parts intersect in at most one point, and therefore their symmetric difference has size at least $2w-2\geq2(w-2)+1.$ Thus Property $\textbf{C}(w-2)$ holds. Property $\textbf{C}(w-2)$ implies that the parts $\mathcal P_{i,j}$ are pairwise disjoint.
\end{proof}

\begin{theorem}\label{cor:rbibd-2}
Let $w\geq4$ and $t\geq\max\{(w-3)\lfloor\frac{w-2}{3}\rfloor-1,w-3\}$. Suppose that $t+1=(w-3)l$ for some positive integer $l$ satisfying $3l\leq w-1$.
If $n>e^{e^{(2w-4)^{12(2w-4)^2}}}$, $n\equiv0\pmod{l(2w-4)}$ and $n\equiv1\pmod{2w-5}$, then
$$C'(n,t,w,w-2)=\frac{n(n-1)}{2((2w-3)(t+1)+l)}.$$
\end{theorem}

\begin{proof}
Apply Theorems~\ref{new-22} and \ref{thm:asy_RBIBD} with $k=2w-4$.
\end{proof}

\section{Concluding remarks}

This paper investigates the upper bounds for $(n,t,w,w-2)$-CPECC codes. Theorem \ref{thm:ub_C'(n,1,w,w-2)} and Theorem \ref{thm:upb_(n,t,w,w-2)-CPECC} examine the cases of $t<w-3$ and $t\geq w-3$, respectively.
Since computing $z_{w-1}(t+1)$ is difficult, Theorem \ref{thm:upb_(n,t,w,w-2)-CPECC_new} provides a new expression of the
upper bound in Theorem \ref{thm:upb_(n,t,w,w-2)-CPECC}. In fact, by a similar proof, ``$t\geq (w-3)\lfloor\frac{w-2}{3}\rfloor-1$" can be replaced by ``$t+1\in\langle w-3,w-2,w-1,w\rangle$" in Theorem \ref{thm:upb_(n,t,w,w-2)-CPECC_new}.
That is, $C'(n,t,w,w-2)\leq\frac{n(n-1)}{2e_{w,t}}$ for $n\geq t+w$, $w\geq4$, $t\geq w-3$, and $t+1\in\langle w-3,w-2,w-1,w\rangle$.
This paper provides some constructions of optimal $(n,t,w,w-2)$-CPECC codes. Theorem \ref{cor-2} applies only when $t\equiv1\pmod2$ and $t<w-3$. An interesting task is to construct optimal $(n,t,w,w-2)$-CPECC codes for $t\equiv0\pmod2$ and $t<w-3$. The construction in Theorem \ref{new-22} requires $t\geq\max\{(w-3)\lfloor\frac{w-2}{3}\rfloor-1,w-3\}$. Note that for $w\geq8$, $(w-3)\lfloor\frac{w-2}{3}\rfloor-1>w-3$. Another interesting task is to construct optimal $(n,t,w,w-2)$-CPECC codes for $w\geq8$ and $w-3\leq t<(w-3)\lfloor\frac{w-2}{3}\rfloor-1$.

\section*{Acknowledgements}

\section*{Declaration on the use of generative AI}

During the preparation of this manuscript, the authors used generative artificial intelligence. The principal ideas and proof outline had already been developed by the authors before the tool was used. The tool assisted with the proof of Theorem~\ref{thm:e_wt=bs+zds}. The authors reviewed and edited all AI-assisted material, verified the mathematical arguments, and take full responsibility for the content of the manuscript.

\appendix
\section{Appendix}\label{app-A}

\textbf{Proof of Theorem \ref{thm:lb_e(G)_(n,t,w,w-2)} }\
Let $n=|V(G)|$. For any $t$-subset $S$, since $G-S\text{ contains a copy of }K_w$, $G-S$ has at least $w$ vertices. Then $n\geq w+t.$

We first deal with the $K_{w+t}$-free setting, and then the unrestricted case.

Suppose that $G$ is $K_{w+t}$-free. We claim that $n\geq 2w+t-1$.  Indeed, assume for contradiction that $n\leq2w+t-2$. Set $H=\overline G$ and $q=n-w-t$.
By $n\geq w+t$, $q\geq0$. Since $G$ is $K_{w+t}$-free, $\omega(G)\leq w+t-1$. Using \eqref{eq:complement-alpha} and \eqref{eq:cover-alpha}, $\tau_v(H)=n-\alpha(H)=n-\omega(G)\geq n-(w+t-1)=q+1.$ Apply Lemma \ref{lem:critical-cover} to $H$ with $r=q+1$. There exists a subgraph $J\subseteq H$ such that
$$\tau_v(J)=q+1 \ \text{ and }\ |V(J)|\leq2(q+1).$$
Note that $n-|V(J)|\geq n-2(q+1)=n-2(n-w-t+1)=2w+2t-2-n\geq t.$
Thus there exists a $t$-set $S\subseteq V(G)\setminus V(J)$. Since $S$ is disjoint from $V(J)$, $J$ is a subgraph of $H-S$. By \eqref{eq:cover-monotone}, $$\tau_v(H-S)\geq\tau_v(J)=q+1.$$
Note that the graph $H-S$ is the complement of $G-S$ on the common vertex set $V(G)\setminus S$. Then $\alpha(H-S)=\omega(G-S)\geq w$ since $G-S\text{ contains a copy of }K_w$. Since $H-S$ has $n-t$ vertices, by \eqref{eq:cover-alpha}, $$\tau_v(H-S)=(n-t)-\alpha(H-S)\leq(n-t)-w=q,$$
a contradiction. Therefore $n\geq 2w+t-1$.

Since $\delta(G)\geq w-1$, apply Lemma \ref{lem:ordering} to $G$ with $k=w-1$. There is a set $X\subseteq V(G)$ such that $G[X]$ is $(w-1)$-colorable and
\begin{align}\label{eq:X-lower}
|X|\geq(w-1)\sum_{v\in V(G)}\frac{1}{d_G(v)+1}.
\end{align}
Note that a $(w-1)$-colorable graph contains no $K_w$. Hence every copy of $K_w$ in $G$ meets $V(G)\setminus X$; in other words, $V(G)\setminus X$ is a $K_w$-transversal. Since $G-S\text{ contains a copy of }K_w$ for any $t$-subset $S$, every $K_w$-transversal has size at least $t+1$, and so $n-|X|=|V(G)\setminus X|\geq t+1$. Thus $|X|\leq n-t-1$. Combining this with \eqref{eq:X-lower} gives
\begin{align}\label{eq:reciprocal-upper}
(w-1)\sum_{v\in V(G)}\frac{1}{d_G(v)+1}\leq n-t-1.
\end{align}

Let $V(G)=\{v_1,v_2,\dots,v_n\}$. Applying the Cauchy--Schwarz inequality to $(\frac{1}{\sqrt{d_G(v_1)+1}},\ldots,\frac{1}{\sqrt{d_G(v_n)+1}})$ and $(\sqrt{d_G(v_1)+1},\dots,\sqrt{d_G(v_n)+1})$,
\begin{align*}
n^2=(1+\cdots+1)^2 \leq (\sum_{i=1}^{n}\frac{1}{d_G(v_i)+1}) (\sum_{i=1}^{n}d_G(v_i)+n).
\end{align*}
Since $\sum_v d_G(v)=2e(G)$,
\begin{align*}
\sum_{v\in V(G)}\frac{1}{d_G(v)+1}\geq\frac{n^2}{2e(G)+n}.
\end{align*}
Together with \eqref{eq:reciprocal-upper}, this yields $\frac{(w-1)n^2}{2e(G)+n}\leq n-t-1$. Since $2e(G)+n>0$ and $n-t-1>0$, $2e(G)+n\geq\frac{(w-1)n^2}{n-t-1}$, and so
\begin{align}\label{eq:Phi-bound}
e(G)\geq\frac n2\left(\frac{n(w-1)}{n-t-1}-1\right).
\end{align}

Let $n_0= 2w+t-1$ and $\Phi(x)=\frac x2\left(\frac{x(w-1)}{x-t-1}-1\right)$ for real $x\geq n_0$. Since $\frac{x^2}{x-t-1}=x+(t+1)+\frac{(t+1)^2}{x-t-1}$,
$$\Phi'(x)=\frac12\left(w-2-\frac{(w-1)(t+1)^2}{(x-t-1)^2}\right).$$
Note that
$$2\Phi'(x)\geq w-2-\frac{(w-1)(w-3)^2}{4(w-1)^2}=w-2-\frac{(w-3)^2}{4(w-1)}=\frac{3w^2-6w-1}{4(w-1)}>0,$$
where the first inequality holds since $x-t-1\geq2(w-1)$ and $t+1\leq w-3$, and the final inequality holds because $3w^2-6w-1>0$ by $w\geq t+4\geq5$.
Hence $\Phi$ is strictly increasing on $[n_0,\infty)$. By \eqref{eq:Phi-bound}, $e(G)\geq\Phi(n)\geq\Phi(n_0)$. Since $e(G)$ is an integer,
\begin{align}\label{eq:rounded-Phi}
e(G)\geq\left\lceil\Phi(n_0)\right\rceil =\left\lceil\frac{(2w+t-1)(2w+t-3)}4\right\rceil=\left\lceil\frac{n_0(n_0-2)}4\right\rceil.
\end{align}
We claim that
\begin{align}\label{eq:ceiling-identity}
\left\lceil\frac{n_0(n_0-2)}4\right\rceil =\binom{\lfloor n_0/2\rfloor}{2}+\binom{\lceil n_0/2\rceil}{2}.
\end{align}
Indeed, for $n_0=2m$, both sides equal $m(m-1)$. If $n_0=2m+1$, then the left-hand side in \eqref{eq:ceiling-identity} is $\left\lceil m^2-\frac14\right\rceil=m^2$, and the right-hand side in \eqref{eq:ceiling-identity} is $\binom m2+\binom{m+1}2=m^2$. Since $\lfloor\frac {n_0}{2}\rfloor=w+\lfloor\frac{t-1}{2}\rfloor$ and $\lceil\frac {n_0}{2}\rceil=w+\lceil\frac{t-1}{2}\rceil$, by \eqref{eq:rounded-Phi} and \eqref{eq:ceiling-identity},
$$e(G)\geq B(w,t)$$
whenever $G$ is $K_{w+t}$-free.

Now we proceed to prove \eqref{eq:ordinary-bound}. Note that
\begin{align*}
\begin{aligned}
4\left(\Phi(n_0)-\binom{w+t}{2}\right)
&=2w^2-t^2-6w-2t+3\\
&\geq2(t+4)^2-t^2-6(t+4)-2t+3=t^2+8t+11>0,
\end{aligned}
\end{align*}
where the first inequality holds since $2w^2-6w$ is increasing for $w\geq t+4\geq5$. Since $B(w,t)=\lceil\Phi(n_0)\rceil$,
\begin{align}\label{eq:B-strict}
B(w,t)>\binom{w+t}{2}.
\end{align}
Hence if $G$ is $K_{w+t}$-free, then $e(G)>\binom{w+t}{2}$. If $G$ contains a copy of $K_{w+t}$, then $e(G)\geq e(K_{w+t})=\binom{w+t}{2}$. To sum up, \eqref{eq:ordinary-bound} holds.

Suppose that equality holds in \eqref{eq:ordinary-bound}. The $K_{w+t}$-free case is impossible by \eqref{eq:B-strict}. Thus $G$ contains a clique $Q\cong K_{w+t}$. Since all $\binom{w+t}{2}$ edges of $G$ are already accounted for by the edges inside $Q$, there are no edges having an endpoint outside $Q$. Any additional vertex would consequently have degree zero, contradicting $\delta(G)\geq w-1>0$. Hence $V(G)=V(Q)$ and $G\cong K_{w+t}$.

Finally, we prove that \eqref{eq:strong-bound} is sharp. Let $a=w+\left\lfloor\frac{t-1}{2}\right\rfloor$, $b=w+\left\lceil\frac{t-1}{2}\right\rceil$, and $G_\star=K_a\mathbin{\dot\cup}K_b$ be the disjoint union of the two complete graphs. Both $a$ and $b$ are at least $w$, so $\delta(G_\star)=\min\{a-1,b-1\}\geq w-1$.
Since $\max\{a,b\}=b\leq w+t-1$, neither component contains $K_{w+t}$, and so $G_\star$ is $K_{w+t}$-free. Every clique in $G_\star$ lies entirely in one of its two components. To destroy every $K_w$ inside a clique $K_m$, one must leave at most $w-1$ vertices and therefore delete at least $m-w+1$ vertices. Consequently, destroying all copies of $K_w$ in both components of $G_\star$ requires at least $(a-w+1)+(b-w+1)=t+1$ deletions. Thus deleting any $t$ vertices from $G_\star$ leaves a copy of $K_w$. Since $e(G_\star)=\binom a2+\binom b2=B(w,t)$, the desired conclusion is obtained. \qed

\section{Appendix}\label{app-B}

\textbf{Proof of Theorem \ref{thm:shadow} }\
We use induction on $n=|V(Q)|$, where $V(Q)$ consists only of vertices lying in at least one hyperedge of $H$.

If $H$ is empty, then $Q$ is empty, $m=0$, $e(Q)=0$ and $z_d(0)=0$, and so the result is immediate. Henceforth suppose that $H$ is nonempty, so $n\geq d+1$.

\medskip
\noindent\textbf{Case 1: $Q$ is disconnected.}

Let the connected components of $Q$ be $Q_1,\dots,Q_k$ with $k\geq2$. Since any two vertices of a hyperedge are adjacent in $Q$, every hyperedge of $H$ lies entirely in one component of $Q$. Let $H_i$ be the subhypergraph consisting of all hyperedges contained in $V(Q_i)$. Then $Q_i$ is a 2-shadow of $H_i$. The vertex sets of different components are disjoint, so the transversal problems in distinct components are independent. Set $m_i=\tau(H_i)$. Then
$$ m=\sum_{i=1}^k m_i\ \text{ and }\ e(Q)=\sum_{i=1}^k e(Q_i). $$
Each component has fewer than $n$ vertices, so the induction hypothesis gives
$$ e(Q_i)\geq bm_i+z_d(m_i)\ \text{ for any }\ 1\leq i\leq k. $$
Summing and applying the subadditivity in \cref{lem:z-properties} give
\begin{align*}
e(Q)\geq b\sum_i m_i+\sum_i z_d(m_i)\geq bm+z_d\left(\sum_i m_i\right)=bm+z_d(m).
\end{align*}

\medskip
\noindent\textbf{Case 2: $Q$ is a clique ($Q$ is connected).}

Write $Q=K_n$. Choose any $d$ vertices to keep and delete the remaining $n-d$ vertices. Since only $d$ vertices remain, no $(d+1)$-element hyperedge can survive. Hence the $n-d$ deleted vertices form a transversal of $H$, and therefore $n\geq d+m$. Since $Q$ is a clique,
\begin{align}
e(Q)=\binom n2\geq\binom{d+m}{2}=bm+\frac{(m-d)(m-d+1)}2.\notag
\end{align}
Set $j=\dist(m,I_1)=\dist\bigl(m,[d-1,d]_{\Z}\bigr)$. If $m\leq d-1$, then $j=d-1-m$; if $m\geq d$, then $j=m-d$. In either case,
$$\frac{(m-d)(m-d+1)}2=\frac{j(j+1)}2\geq j.$$
Since $I_1$ is one of the intervals over which the minimum defining $z_d(m)$ is taken, $z_d(m)\leq j$. Thus
$$e(Q)\geq bm+z_d(m).$$

\medskip
\noindent\textbf{Case 3: $Q$ is connected and $Q$ is not a clique, and $\Delta(Q)\leq b$.}

We first prove that $Q$ has a proper $b$-colouring. If $\Delta(Q)\leq b-1$, then the standard greedy colouring algorithm uses at most $\Delta(Q)+1\leq b$ colours. For $\Delta(Q)=b$, $Q$ is not an odd cycle because $\delta(Q)\geq d\geq3$. Since $Q$ is not a clique, $\chi(Q)\leq b$ by Brooks' theorem. Thus $Q$ admits a proper colouring with at most $b$ colours. Add empty colour classes if fewer than $b$ colours are actually used, and write the colour-class sizes in nonincreasing order as $x_1\geq x_2\geq\cdots\geq x_b\geq0$.

The union of the first $d$ colour classes induces a graph that is properly $d$-colourable, hence contains no $K_{d+1}$. Since every hyperedge of $H$ induces a $K_{d+1}$ in $Q$, the first $d$ colour classes cannot contain any hyperedge of $H$. Therefore, the union of the remaining $b-d=d-1$ colour classes is a transversal of $H$, and $$\sum_{i=d+1}^{b}x_i\geq m.$$

Since $x_{d+1}$ is the largest among these last $d-1$ classes, $(d-1)x_{d+1}\geq\sum_{i=d+1}^{b}x_i\geq m$, and so $x_{d+1}\geq\left\lceil\frac{m}{d-1}\right\rceil.$
Since $x_1,\dots,x_d\geq x_{d+1}$,
\begin{equation}\label{eq:n-colour-lower}
n=\sum_{i=1}^{b}x_i\geq d\left\lceil\frac{m}{d-1}\right\rceil+\sum_{i=d+1}^{b}x_i\geq d\left\lceil\frac{m}{d-1}\right\rceil+ m.
\end{equation}

Now, for each $v\in V(Q)$, set $c_v=b-1-d_Q(v)$. Note that $-1\leq c_v\leq d-2$ by $d\leq d_Q(v)\leq b$, and
$$\frac{(d_Q(v)-(2d-2))(d_Q(v)-(2d-1))}{2(d_Q(v)+1)}=\frac{c_v(c_v+1)}{2(b-c_v)}.$$
For a real number $x$, write $x_+=\max\{x,0\}$, and define
$$P=\sum_{v\in V(Q)}(c_v)_+.$$
We claim that for $c_v\geq1$,
\begin{equation}\label{eq:c-ineq}
\frac{c_v(c_v+1)}{2(b-c_v)}\geq\frac{c_v}{b-1}.
\end{equation}
Indeed, after dividing by the positive number $c_v$ and clearing the positive denominators, \eqref{eq:c-ineq} is equivalent to $(c_v+1)(b-1)\geq2(b-c_v),$ and the difference between the left-hand side and the right-hand side is $(c_v-1)(b+1)\geq0.$ Thus $(c_v)_+=c_v\leq (b-1)\frac{(d_Q(v)-(2d-2))(d_Q(v)-(2d-1))}{2(d_Q(v)+1)}$ for $c_v\geq1$.
If $c_v=0$ or $c_v=-1$, then
$$\frac{(d_Q(v)-(2d-2))(d_Q(v)-(2d-1))}{2(d_Q(v)+1)}=\frac{c_v(c_v+1)}{2(b-c_v)}=0=\frac{(c_v)_+}{b-1}.$$
To sum up, $(c_v)_+\leq (b-1)\frac{(d_Q(v)-(2d-2))(d_Q(v)-(2d-1))}{2(d_Q(v)+1)}$. Recall that $D=e(Q)-bm$. Hence 
\begin{equation}\label{eq:P-D}
P\leq(b-1)\sum_{v\in V(Q)}\frac{(d_Q(v)-(2d-2))(d_Q(v)-(2d-1))}{2(d_Q(v)+1)}\leq(b-1)D,
\end{equation}
where the last inequality holds by Lemma \ref{lem:random} (b). On the other hand, since $c_v\leq(c_v)_+$ for every $v$, $\sum_v c_v\leq P$, and so
\begin{equation}\label{eq:degree-lower}
2e(Q)=\sum_v d_Q(v)=(b-1)n-\sum_v c_v\geq(b-1)n-P.
\end{equation}
By using $e(Q)=bm+D$, \eqref{eq:P-D} and \eqref{eq:degree-lower}, we have $2bm+2D\geq(b-1)n-(b-1)D$. That is,
\begin{equation}\label{eq:D-lower-n}
(b+1)D\geq(b-1)n-2bm.
\end{equation}

Write $m$ by Euclidean division by $d-1$:
$$m=q(d-1)+r,\text{ where }q\in\Z_{\geq0}\text{ and }0\leq r<d-1.$$
If $r=0$, then $m=q(d-1)\in I_q$, and so $z_d(m)=0$. Lemma \ref{lem:random} (b) implies $D\geq z_d(m)$. Suppose now that $r>0$. Then
$$\left\lceil\frac{m}{d-1}\right\rceil=q+1.$$
By \eqref{eq:n-colour-lower}, we have
\begin{align*}
n\geq m+d(q+1)=q(d-1)+r+dq+d=(2d-1)q+d+r=bq+d+r.
\end{align*}
Substituting this into \eqref{eq:D-lower-n} and using $b-1=2(d-1)$ gives
\begin{align*}
(b+1)D
&\geq(b-1)(bq+d+r)-2b\bigl(q(d-1)+r\bigr)\\
&=(b-1)(d+r)-2br=2d(d-1)-(b+1)r=(b+1)(d-1-r).
\end{align*}
Therefore
\begin{equation}\label{eq:D-r}
D\geq d-1-r.
\end{equation}
The left endpoint of $I_{q+1}$ is
$$(q+1)(d-1)=q(d-1)+r+(d-1-r)=m+(d-1-r).$$
Hence $z_d(m)\leq d-1-r$. Together with \eqref{eq:D-r}, this yields $D\geq z_d(m)$, as required.

\medskip
\noindent\textbf{Case 4: $Q$ is connected and $\Delta(Q)\geq b+1$.}

Choose a vertex $v\in V(Q)$ of maximum degree. Let $H'=H-v$ and $m'=\tau(H')$, and let $Q'$ be the $2$-shadow of $H'$. If $T$ is a transversal of $H$, then $T\setminus\{v\}$ is a transversal of $H'$, and so $m'\leq m$. Conversely, if $T'$ is a minimum transversal of $H'$, then $T'\cup\{v\}$ is a transversal of $H$, and so $m\leq m'+1$. To sum up,
\begin{equation}
m'\in\{m-1,m\}.\notag
\end{equation}

Since every edge of $Q$ incident with $v$ disappears from $Q'$, we have
\begin{equation}
e(Q)-e(Q')\geq d_Q(v)\geq b+1.\notag
\end{equation}
Because the induction hypothesis applies to $H'$, we have $e(Q')\geq bm'+z_d(m')$. If $m'=m$, then
\begin{align*}
e(Q)\geq e(Q')+b+1\geq bm+z_d(m)+b+1\geq bm+z_d(m).
\end{align*}
If $m'=m-1$, then
\begin{align*}
e(Q)\geq e(Q')+b+1&\geq b(m-1)+z_d(m-1)+b+1\\
&=bm+z_d(m-1)+1\geq bm+z_d(m),
\end{align*}
where the last inequality follows from \cref{lem:z-properties}.

The four cases exhaust all possibilities for $Q$, so this completes the proof.
\qed

\end{document}